\documentclass[a4paper,11pt,twoside]{article}
\usepackage{pst-fill,pst-grad}
\usepackage{textcomp}
\usepackage[english]{babel} 
\usepackage[utf8]{inputenc} 
\usepackage{graphicx}
\usepackage{amsmath} 
\usepackage{float}  
\usepackage{fancyhdr} 
\usepackage[matrix,arrow,curve]{xy}
\usepackage{pstricks} 
\usepackage{amsmath,amsfonts,verbatim,afterpage,theorem,euscript,mathrsfs,amssymb}
\usepackage{amsfonts}
\usepackage{amssymb}
\usepackage{array} 
\usepackage{dsfont}
\usepackage[colorlinks=true,linkcolor=blue,citecolor=red]{hyperref}
\usepackage{authblk}
\usepackage{color}
\usepackage{nicefrac}
\usepackage[titletoc]{appendix}
\usepackage{titlesec}
\newtheorem{Proposition}{Proposition}[section]

\newtheorem{Lemme}{Lemma}[section]

\newtheorem{TheoremeP}{Theorem}
\newtheorem{CorollaireP}{Corollary}

\newcommand{\R}{\mathbb{R}}
\newcommand{\N}{\mathbb{N}}

\newcommand{\norm}[1]{\left \Vert #1 \right \Vert}
\newcommand{\abs}[1]{\left \vert #1\right \vert}
\newcommand{\func}[3]{#1\colon  #2  \to  #3}
\def \R{\mathbb{R}}
\def \Rn{\mathbb{R}^n}

\def \finpv{\hfill $\blacksquare$  }
\def \pv{{\bf{Proof.}}~} 

\def \ds{\displaystyle}

\title{\bf  From anomalous to classical diffusion in a non-linear heat equation}

\author[1]{ Oscar Jarr\'in\footnote{corresponding author: oscar.jarrin@udla.edu.ec}}
\author[2]{Geremy Loacham\'in\footnote{geremy.loachamin@uni.lu}}

\affil[1]{\scriptsize Escuela de Ciencias Físicas y Matemáticas, Universidad de Las Américas, Vía a Nayón, C.P.170124, Quito, Ecuador.}
\affil[2]{\scriptsize Faculty of Science,  University of Luxembourg, Maison du Nombre, 6 Avenue de la Fonte, Esch-sur-Alzette L-4364, Luxembourg.}
\date{\today}   
\begin{document}
\maketitle	 
%%%%%%%%%%%%%%%%%%%%%%%%%%%%%%%%%%%%%%%%%%%%%%
\begin{abstract}
In this paper, we consider the heat equation with the natural polynomial non-linear term; and  with two  different cases in the diffusion term. The first case (anomalous diffusion) concerns to the fractional Laplacian operator  with  parameter $1<\alpha <2$, while  the second case (classical diffusion) involves the classical Laplacian operator. When $\alpha \to 2$,  we prove  the uniform convergence of the solutions of the anomalous diffusion case to a solution of the classical diffusion case. Moreover, we rigorously derive a convergence rate, which was experimentally exhibited in previous related works. \\[3mm]
\end{abstract}
%\tableofcontents
%%%%%%%%%%%%%%%%%%%%%%%%%%%%%%%%%%%%%%%%%%%%%%
%%%%%%%%%%%%%%%%%%%%%%%%%%%%%%%%%%%%%%%%%%%%%%
\section{Introduction}\label{Sec:Intro} 
In this paper, we consider the following multi-dimensional, nonlinear  and anomalous diffusion heat equation in the whole space $\Rn$ with $n \geq 1$:
\begin{equation}\label{Eq}
\partial_t u  + (-\Delta)^{\alpha /2}\, u + \eta \cdot  \nabla (u^b)=0, \qquad 1<  \alpha  < 2, \quad b \in \mathbb{N} \quad \mbox{with}\quad  b \geq 2.
\end{equation}
Here, the function  $u:[0,+\infty[\times \Rn\to \R$ is  the solution, and  $(-\Delta)^{\alpha /2}\, u$ is the \emph{anomalous} diffusion term which is given by the fractional Laplacian operator $(-\Delta)^{\alpha/2}$. We recall that this operator is  defined  in the Fourier level  by $\ds{\widehat{(-\Delta)^{\alpha/2} u}(t,\xi)= c_{n,\alpha}\, \vert \xi \vert^{\alpha}\widehat{u}(t,\xi)}$. Moreover,  in the spatial variable, the fractional Laplacian operator  is defined as the following non-local operator:   
$$ (-\Delta)^{\alpha /2} u(t,x)= c_{n,\alpha}\ \text{\textbf{p.v.}}\int_{\Rn} \frac{u(t,x)-u(t,y)}{\vert x-y \vert^{n+\alpha}} dy,$$
 where $\textbf{p.v.}$ denotes the principal value. Finally,  $\eta \in \Rn$ is a fixed vector, and moreover,  the parameter  $b \in \mathbb{N}$  in the nonlinear term verifies    $b\geq 2$. 
 
 \medskip 
 
 We may observe that this \emph{highly nonlinear term} essentially behaves as the derivative of a polynomial of degree $b$ in the variable $u$. Thus, this term  agrees with the classical assumption for the non-linearity in the qualitative study of the heat equation. See, for instance,     \cite{Biler1,BiTaWo,BiKarWo,BrandKarch,IgnatRossi} and the references therein.    
 
\medskip  

Nonlinear  evolution  PDEs involving the fractional Laplacian, which describe the \emph{anomalous} or $\alpha-$Lévy stable diffusion,  have been extensively studied in the physical and mathematical points of view. From the physical point of view, and for $b=2$, the equation (\ref{Eq}) deals with a generalized Burgers-type equation \cite{BiTaWo} which has been largely  used to model a variety of physical phenomena such as, for example,  the anomalous homogeneous turbulence \cite{FuSuWo}, applications to hydrodynamics and  statistical  mechanics \cite{ShZaFri},  and moreover, applications to molecular biology in the modeling of  growth of molecular interfaces \cite{Zas2}. In the latter application, the general algebraic non-linear term $u^b$, with $b \in \mathbb{N}^*$  and $b\geq 2$,  provides a good model for multi-particle interactions. For more references, see the book \cite{SaiWo2}.

\medskip

From the mathematical point of view, when the solution $u(t,\cdot)$ is considered as  the density of a probability distribution for every $t > 0$, the equation (\ref{Eq}) has an important probabilistic interpretation in the theory of nonlinear Markov processes and propagation of chaos. See, \emph{e.g.}, the works \cite{FuWo2}, \cite{Jourdan} and the references therein.

\medskip

Getting back to the expression (\ref{Eq}), we observe that for each value of the parameter $1<\alpha <2$ in the fractional Laplacian operator $(-\Delta)^{\alpha /2}$ we get a corresponding  fractional PDE.  We thus  denote   by $u_\alpha(t,x)$ the corresponding solution of each  equation and 
the main objective of this paper is to understand the asymptotic behavior of the family of functions $u_\alpha(t,x)$ when the parameter $\alpha$ goes to $2$.  This question was pointed out from the experimental point of view in  \cite{FuSuWo,Zas2} and has some interesting applications in these physical and biological models. Our aim is then to provide a rigorous mathematical framework  to give an answer.  

\medskip

Formally, we may observe that  if in the expression (\ref{Eq}) we set $\alpha=2$, then we  get a \emph{classical   diffusion equation}  involving the Laplacian operator: 
\begin{equation}\label{Eq2}
\partial_t u  - \Delta \, u + \eta \cdot  \nabla (u^b)=0.
\end{equation} 
Consequently, if $u(t,x)$ denotes a  solution of the equation above, we are  interested in providing a rigorous understanding of the expected convergence $u_\alpha(t,x) \to u(t,x)$, when $\alpha \to 2$. It is worth mentioning  although this problem is easily formulated, it is not a trivial study since for each value of the parameter $\alpha$ we have a different fractional PDE depending on this parameter.   

\medskip

In the  particular   case  of the following  \emph{linear equation} in   a smooth and bounded domain $\Omega\subset \Rn$:
\begin{equation}\label{EqLin}
\partial_t u_\alpha+ (-\Delta)^{\alpha /2}\, u_\alpha =f_\alpha, \quad 0< \alpha  < 2,
\end{equation}  
and where the function $f_\alpha(t,x)$ does not depend on the solution $u_\alpha$,   this convergence  problem  was   studied by U. Biccari \& V. Hern\'andez-Santamar\'ia  in   \cite{Biccari}. For a time $0<T<+\infty$, the authors consider a family of functions   $f_\alpha \in L^2(0,T, H^{-\alpha}(\Omega))$, which is uniformly bounded with respect to the parameter $\alpha$:  $\Vert f_\alpha(t,\cdot)\Vert_{H^{-\alpha} (\Omega)} \leq C $, and such  that when $\alpha \to 2$ we have the  convergence $f_\alpha(t,\cdot)\to f(t,\cdot)$ in  the weak topology of the space $H^{-1}(\Omega)$.  Then, by using  a compactness argument (due to the boundness of the domain $\Omega$) it is shown that the \emph{weak solutions} of equation (\ref{EqLin}) converge in the strong topology of the space $L^2(0,T, H^{1-\delta}_{0}(\Omega))$ (with $0<\delta \leq 1$) to a \emph{weak solution} of the corresponding linear heat equation with datum $f$. Moreover,  when  $\delta=1$,  in the setting of the space $L^2(0,T, L^{2}(\Omega))$  the authors  \emph{numerically} obtain a convergence rate of the order $\vert 2-\alpha \vert^{1/2}$.

\medskip

On the other hand,  L. Ignat \& J.D. Rossi proved in \cite{IgnatRossi}, among other things, that \emph{weak solutions}  $u(t,x)$ fo the nonlinear heat equation (\ref{Eq2}) can be obtained as the limit   (when $\varepsilon\to 0^{+}$) of the \emph{weak solutions}  to the following nonlocal convection-diffusion  equation in the whole space $\Rn$:
\begin{equation}\label{EqJ}
\partial_{t} u_\varepsilon  + \frac{1}{\varepsilon^2}(J_\varepsilon\ast u_\varepsilon -u_\varepsilon)+ \frac{1}{\varepsilon} (G_\varepsilon \ast  u^{b}_{\varepsilon} -u^{b}_{\varepsilon} )=0, \quad \varepsilon>0.  
\end{equation} This equation has the same scaling properties of the equation  (\ref{Eq2}) and  here, for suitable non-negative functions $J\in \mathcal{S}(\Rn)$  and $G \in \mathcal{S}(\Rn)$,  we have $J_\varepsilon(x)= \frac{1}{\varepsilon^n} J (x / \varepsilon) $ and $G_\varepsilon(x)= \frac{1}{\varepsilon^n} G (x / \varepsilon)$ respectively. Moreover,  $J$  is a radially symmetric function and the key assumption is that its Fourier transform  $\widehat{J}(\xi)$ satisfies  the following condition: 
\begin{equation}\label{Cond}
  \frac{1}{2}\partial^{2}_{\xi_i} \widehat{J_\varepsilon}(0)=1, \qquad  i=1,\cdots, n,    
\end{equation}
which is similarly satisfied for the symbol $\vert \xi \vert^2$ of the classical Laplacian operator. In this setting, by using some sharp estimates of the kernel associated to the linear problem, and moreover, by setting  the  vector $\eta=(\eta_1, \cdots, \eta_n)$ in the equation (\ref{Eq2}) as $\ds{\eta_i= \int_{\Rn} x_i G(x)dx}$, for all time $0<T<+\infty$ it is proven  the following convergence result in the natural framework (due to the Plancherel's identity) of the Lebesgue space $L^2(\Rn)$:
\begin{equation}\label{Conv}
 \lim_{\varepsilon \to 0^{+}} \sup_{0\leq t \leq T} \Vert u_\varepsilon(t,\cdot)-u(t,\cdot)\Vert_{L^2(\Rn)}=0.    
\end{equation}

Nevertheless, we remark that  this result cannot  be applied to the case of the equation (\ref{Eq}). Indeed, since the symbol $\vert \xi \vert^\alpha$ of the operator  $(-\Delta)^{\alpha / 2}$  does not verify the key condition (\ref{Cond}) the  nonlocal   diffusion operator $1/\varepsilon^2(J_\varepsilon \ast (\cdot) - I_d)$ (where $I_d$ is the identity operator) \emph{does not contain} the fractional Laplacian operator $(-\Delta)^{\alpha / 2}$  as a particular case. Moreover, one  the main property of the approximated equation  (\ref{EqJ}) is the same scaling of the equation (\ref{Eq2}), which is not verified when considering the fractional Laplacian operator.

\medskip 

In this work, we  will use a  different approach. For any time $0<T<+\infty$, in the framework the space $L^{\infty}([0,T]\times \R^n)$  we shall study the convergence (in the strong topology) of the \emph{strong (mild) solutions} $u_{\alpha}(t,x)$ for  the anomalous diffusion equation (\ref{Eq}) (given in the expression (\ref{u-alpha})) to a \emph{strong solution} $u(t,x)$ for the classical diffusion equation (\ref{Eq2}) (given in the expression (\ref{u-2})). See our main result given in Theorem \ref{Th-Unif-Conv} for the details. This uniform convergence  also  allows us to prove a strong convergence in the $L^{p}_{t} L^{q}_{x}$ spaces (see the Corollary \ref{Col-1}).  

\medskip

 Our method is based on two key ideas. On the one hand,  we   study  the convergence of the fundamental solution $p_\alpha(t,x) $ associated with the  fractional linear heat equation  (see the expression (\ref{pa}) for a definition)   to the  heat kernel $h(t,x)$. On the other hand, we prove some uniform estimates with respect to the parameter $\alpha$ for the family of functions $u_{\alpha}(t,x)$.  

\medskip

Finally, we think that in a further research  our method could be adapted to the  case when the fractional Laplacian operator $(-\Delta)^{\alpha/2}$ in the equation (\ref{Eq}) is substituted by a more general Lévy-type operator $\mathcal{L}^\alpha$. For a definition and some well-known properties of this latter operator, we refer to the book \cite{Jacob}.   

\section{The main result} 

Let us consider the Cauchy problem for both  \emph{anomalous} (when $1<\alpha <2$) and \emph{classical} (when $\alpha=2$) nonlinear heat equation: 
\begin{equation}\label{Eq3}
\left\{ \begin{array}{ll}\vspace{2mm}
\partial_t u_\alpha + (-\Delta)^{\alpha /2}\, u_\alpha+ \eta \cdot  \nabla (u^{b}_{\alpha})=0, \quad 1< \alpha \leq 2,\\
u_\alpha(0, \cdot)=u_{0,\alpha}. 
\end{array} \right.
\end{equation}
Well-posedness (WP) issues  for this equation have been studied in several works \cite{Biler1,Dro,Dro2} and it is well-known that for an initial datum $u_{0} \in L^1(\Rn)$  the initial value problem (\ref{Eq3}) has a unique solution $u_\alpha \in \mathcal{C}([0,+\infty[, L^1(\Rn))$ which verifies 
\begin{equation}\label{Bound-L1}
\Vert u_\alpha (t,\cdot) \Vert_{L^1} \leq \Vert u_{0} \Vert_{L^1}.     
\end{equation}
Moreover, for $1\leq p \leq +\infty$ this solution also verifies $\ds{u_\alpha\in \mathcal{C}(]0,+\infty[, W^{1,p}(\Rn))}$, and the following estimate holds:
\[   \Vert u_\alpha(t,\cdot) \Vert_{L^p} \leq C\, t^{-\frac{n}{\alpha}(1-\frac{1}{p})}\, \Vert u_{0}\Vert_{L^1}.\]
Finally, under the additional assumption on the initial datum: $u_{0} \in L^1 \cap  L^p(\Rn)$ the corresponding solution verifies $\ds{u_\alpha \in \mathcal{C}([0,+\infty[, L^p(\Rn))}$, and for all time $t\geq 0$ we have the estimate 
\[ \Vert u_\alpha (t,\cdot)\Vert_{L^p} \leq \Vert u_{0} \Vert_{L^p}.\]

As mentioned, our aim is to study the convergence (when $\alpha \to 2$) of mild solutions of the equation (\ref{Eq3}). For technical  reasons, principally due to the study of the limit concerning  the highly nonlinear term $\nabla (u^{b}_{\alpha}) \to \nabla (u^{b}_{2}) $ (recall that we assume $b$  an integer such that $b \geq 2$),  we shall need more regularity than the one given by the space $W^{1,p}(\R^n)$.  For this,  we shall consider  initial data belonging to  the space $L^1 (\R^n) \cap H^s(\R^n)$ with $s > n / 2$.  

\medskip

The global-well posedness in the setting of the space $L^1 (\R^n) \cap H^s(\R^n)$ is rather standard but, to our knowledge, this fact has not been proven before. Consequently, only  for the   completeness of this paper,  we state and we will  give  a proof of the following theorem. We emphasize that the only novelty is the gain of  regularity   given in the expression (\ref{Regularity-solutions}) below, which follows from the additional hypothesis $u_0 \in H^s(\R^n)$.  

\begin{TheoremeP}\label{Th-WP} Let $1<\alpha \leq 2$.  For $s>n/2$,  let 
$u_{0}  \in L^1 \cap  H^s(\Rn) $   be an initial datum. Then there exists a unique mild solution 
\[ u_\alpha \in \mathcal{C}([0,+\infty[, L^1 \cap H^s(\Rn)), \]
of the equation (\ref{Eq3}). Moreover,  this solution is regular: 
\begin{equation}\label{Regularity-solutions}
u_\alpha \in \mathcal{C}^{1}(]0, +\infty[, \mathcal{C}^{\infty}(\Rn)), 
\end{equation}
and it verifies the equation (\ref{Eq3}) in the classical sense. 
\end{TheoremeP}  

\medskip 

We study now the convergence of mild solutions for the equation (\ref{Eq3}) when $\alpha \to 2^{-}$.  For the fractional case (when $\alpha <2$) we consider  $(u_{0,\alpha})_{1<\alpha < 2} \subset L^1\cap H^s(\Rn)$  family of initial data;  and we shall denote by  $u_\alpha(t,x)$ the corresponding arising solution  given by Theorem \ref{Th-WP}. Moreover, for the classical case (when $\alpha=2$) we similarly consider an initial datum and its corresponding solution $u_{0,2} \in L^1\cap H^s(\Rn)$ and  $u_2(t,x)$ respectively.

\medskip

We shall assume the  following strong  convergence on the initial data: 
\begin{equation}\label{Conv-Data}
u_{0,\alpha} \to u_{0,2}, \quad \alpha \to 2^{-}, \qquad \mbox{in} \quad  L^1 \cap H^s(\R^n).  
\end{equation}
 By the Sobolev embedding  $H^s(\Rn)\subset L^{\infty}(\Rn)$ (since $s>n/2$)  we also have $u_{\alpha,0} \in L^{\infty}(\R^n)$,  $u_{0,2}\in L^{\infty}(\R^n)$ and the convergence above also holds true in the space $L^{\infty}(\R^n)$. Thus, for the corresponding family of solutions $(u_\alpha)_{1<\alpha \leq 2}$  we will  study  the \emph{uniform convergence}: 
\begin{equation}\label{Conv-Sol}
 u_\alpha(t,x) \to u_2(t,x), \quad \alpha \to 2^{-}, \qquad \mbox{in} \quad L^{\infty}([0,T]\times \Rn) \quad \mbox{for any} \quad 0<T<+\infty. 
 \end{equation}
Moreover, we are also interested in studying the convergence rate in (\ref{Conv-Sol}). For this, we introduce a parameter $\gamma>0$ and we shall assume the estimate (\ref{Conv-Rate-Data} below, which is a given convergence rate of the initial data in the space $L^\infty(\R^n)$. Our aim is then  to study when the family of solutions follows this prescribed convergence rate.   In this setting, our main result reads as follows:  

\begin{TheoremeP}\label{Th-Unif-Conv} Let $(u_{0,\alpha})_{1<\alpha \leq 2}$ be a family of initial data such that for all   $1<\alpha \leq 2$ we have 
$u_{0,\alpha} \in L^1 \cap H^s(\Rn)$.  Let $ (u_\alpha)_{1<\alpha \leq  2}$ be corresponding family of solutions to  the equation (\ref{Eq3}) given by  Theorem \ref{Th-WP}.  

\medskip

We assume the convergence given in (\ref{Conv-Data}), and moreover, for a parameter $\gamma>0$   we assume   the estimate
\begin{equation}\label{Conv-Rate-Data}
 \Vert u_{0,\alpha}-u_{0,2}\Vert_{L^\infty} \leq {\bf c}\, (2-\alpha)^{\gamma},    
\end{equation}
where ${\bf c}>0$ is  a given  generic constant. Then, there exists $0<\varepsilon \ll 1$,  and there exists a   constant ${\bf C}>0$, which depends on   the parameters $\eta$ and $b$ in the equation (\ref{Eq3}), the initial $u_{0,2}$, the quantity $\varepsilon$   and the constant ${\bf c}$, such that for all  $1+\varepsilon < \alpha< 2$ the following estimate holds: 
\begin{equation}\label{Conv-Rate-Sol}
\sup_{0\leq t \leq T} \Vert  u_\alpha (t,\cdot)-  u_{2}(t,\cdot)\Vert_{L^{\infty}} \leq {\bf C}\, (1+T+T^2)\, \max\Big((2-\alpha)^\gamma, 2-\alpha \Big), \quad 0<T<+\infty. 
\end{equation}		 
\end{TheoremeP}	

\medskip

Some remarks are in order here. First note that  our approach allows us obtain  a uniform convergence (in both  the temporal and the spatial variables) which   is not studied in the previous related works  \cite{Biccari} and \cite{IgnatRossi}. Moreover, it is interesting to observe the convergence rate given in the estimate (\ref{Conv-Rate-Sol}) is  determined  by a  \emph{competition} between the   quantities  $(2-\alpha)^\gamma$ and $(2-\alpha)$.

\medskip

In order to make a deeper discussion of this fact, let us briefly explain the general idea of the proof. As pointed out, we shall consider mild solutions of the equation (\ref{Eq3}). Thus,  for $1<\alpha <2$ we have 
\begin{equation}\label{u-alpha}
u_\alpha(t,\cdot)=p_\alpha(t,\cdot) \ast u_{0,\alpha}  + \int_{0}^{t}p_\alpha(t-s,\cdot)\ast \eta \cdot \nabla(u^{b}_{\alpha})(s,\cdot) ds,    
\end{equation}
where the kernel $p_\alpha(t,x)$ is given in (\ref{pa}); while for $\alpha =2$ we  have
\begin{equation}\label{u-2}
u_2(t,\cdot)= h(t,\cdot) \ast u_{0,2} + \int_{0}^{t}h(t-s,\cdot)\ast \eta \cdot \nabla(u^{b}_{2})(s,\cdot) ds,     
\end{equation} where $h(t,x)$ always denotes the heat kernel.   The estimate (\ref{Conv-Rate-Sol}) is then obtained by the following estimates
\begin{equation*}
\left\Vert p_\alpha(t,\cdot) \ast u_{0,\alpha}  - h(t,\cdot) \ast u_{0,2}  \right\Vert_{L^\infty} \lesssim \max\Big((2-\alpha)^\gamma, 2-\alpha \Big),
\end{equation*}
and 
\begin{equation*}
\left\Vert \int_{0}^{t}p_\alpha(t-s,\cdot)\ast \eta \cdot \nabla(u^{b}_{\alpha})(s,\cdot) ds- \int_{0}^{t}h(t-s,\cdot)\ast \eta \cdot \nabla(u^{b}_{2})(s,\cdot) ds  \right\Vert_{L^\infty} \lesssim \max\Big((2-\alpha)^\gamma, 2-\alpha \Big),
\end{equation*}
on the linear and the nonlinear terms respectively.  For the sake of simplicity, we will only explain more in detail the estimates on the linear term. Of course the estimates for the nonlinear term are much more delicate, but they follow some similar ideas.   We split the linear term as 
\begin{equation*}
\left\Vert p_\alpha(t,\cdot) \ast u_{0,\alpha}  - h(t,\cdot) \ast u_{0,2}  \right\Vert_{L^\infty} \leq \left\Vert ( p_\alpha(t,\cdot) - h(t,\cdot))   \ast u_{0,\alpha} \right\Vert_{L^\infty} + \left\Vert h(t,\cdot) \ast (u_{0,\alpha}   -  u_{0,2}) \right\Vert_{L^\infty},
\end{equation*}
where we have 
\begin{equation*}
\left\Vert ( p_\alpha(t,\cdot) - h(t,\cdot))   \ast u_{0,\alpha} \right\Vert_{L^\infty} \lesssim (2-\alpha) \qquad \mbox{and} \qquad  \left\Vert h(t,\cdot) \ast (u_{0,\alpha}   -  u_{0,2}) \right\Vert_{L^\infty} \lesssim (2-\alpha)^{\gamma}. 
\end{equation*}
Here,  the  quantity $(2-\alpha)^\gamma$  is the convergence rate assumed for the initial data, while the quantity $(2-\alpha)$  is   the convergence rate of the kernels  $p_\alpha(t,x) \to h(t,x)$  (when $\alpha \to 2^{-}$), which is proven in  the Lemma \ref{Lema-Tech-2} below. 

\medskip

Since   we have $1 < \alpha <2$ and therefore $0<2-\alpha< 1$, the estimate (\ref{Conv-Rate-Sol}) yields the following conclusions  by considering two cases of the parameter $\gamma$.
\begin{enumerate}
\item[$\bullet$] {\bf  The case $0<\gamma \leq 1$}. Here we have $\ds{\max\Big((2-\alpha)^\gamma, 2-\alpha \Big) = (2-\alpha)^\gamma}$, and consequently, the solutions $u_\alpha(t,x)$ converge to the solution $u_2(t,x)$ with the same convergence rate as that of the initial data. 

\item[$\bullet$] {\bf  The case $\gamma > 1$}. In this case, we have $\ds{\max\Big((2-\alpha)^\gamma , 2-\alpha \Big) = 2-\alpha}$. Then,  it is interesting to observe that the convergence rate of the solutions does not follow the one of initial data. More precisely,  the solutions  $u_\alpha(t,x)$ converge to the solution $u_2(t,x)$ with a rate of order $2-\alpha$, which is \emph{slower} than the convergence rate of the initial data $(2-\alpha)^\gamma$. 
\end{enumerate}
Summarizing, the increasing of the parameter parameter $\gamma$ makes the assumption (\ref{Conv-Rate-Data}) strong but not the result given in (\ref{Conv-Rate-Sol}). This  phenomenological effect is given by the convergence rate  of the kernels   $p_{\alpha}(t,\cdot) \to h(t,\cdot)$. 

\medskip 

On the other hand, as mentioned in the introduction, the convergence result given in Theorem \ref{Th-Unif-Conv} also allows us to study the convergence (\ref{Conv-Sol}) in the following Lebesgue spaces.

\begin{CorollaireP}\label{Col-1} With the same hypothesis  of Theorem \ref{Th-Unif-Conv},  for all $1\leq p \leq +\infty$ and  $1<q<+\infty$  we have the estimate:
\[  \Vert u_\alpha- u_2 \Vert_{L^{p}((0,T], L^q(\Rn))} \leq C_{p,q}\, (1+T+T^2)\,  \max\Big((2-\alpha)^{\gamma\left( 1-\frac{1}{q}\right)}, (2-\alpha)^{1-\frac{1}{q}}\Big), \quad 1+\varepsilon < \alpha < 2.\] 
\end{CorollaireP}

We observe that in the framework of $L^{p}_{t}L^{q}_{x}-$spaces, the convergence rate is only driven by the parameter $q$, which describes the decaying properties of solutions in the spatial variable.  Moreover, by setting the parameter    $\gamma=1$ and with the particular values $p=q=2$, we obtain the following convergence rate: $ \Vert u_\alpha - u_2 \Vert_{L^{2}_{t}L^{2}_{x}} \lesssim  (2-\alpha)^{1/2}$,   which was experimentally obtained in \cite{Biccari} for the particular linear case (when $\eta=0$)  of the equation (\ref{Eq3}). 

\medskip  

To close this section, let us make the following final comments. First note that in this work we have restricted  ourselves in the case when the parameter $\alpha$ verifies $1<\alpha <2$, however, our results  are also valid for the case $\alpha>2$ with minor  technical modifications. 

\medskip

 The lower constraint $1+\varepsilon < \alpha$ (with $0<\varepsilon \ll 1$) given in Theorem \ref{Th-Unif-Conv} is essentially technical, due to estimates involving the expression   $\frac{1}{1-1/\alpha}$ (see, for instance, the estimate  in Proposition \ref{Prop-Non-Lin} below). Consequently, our result left open the convergence problem when $\alpha \to 1^{+}$ which is also interesting and could be a matter of further research. 

\medskip

Finally, we think that our method explained above could be also adapted to study the convergence given in (\ref{Conv-Sol}) within the framework of other functional spaces, provided that we assume some  natural hypothesis on  the initial data.   

\medskip 

\textbf{Organization of the paper.}  In Section \ref{Sec:Known-Facts} we recall some well-known facts on the linear fractional heat equation that we will use in the next sections. Section \ref{Sec:WP-Reg} is devoted to the proof of Theorem \ref{Th-WP}, while in Section \ref{Sec:Anom-to-local}  we give a proof of Theorem \ref{Th-Unif-Conv} and Corollary \ref{Col-1}. 
 
\section{Some well-known facts}\label{Sec:Known-Facts}
In this section, for the completeness of this paper, we  summarize some well-known facts on the linear and homogeneous  fractional heat equation: 
\[ \partial_t p_\alpha + (-\Delta)^{\alpha/2} p_\alpha =0, \quad 1<\alpha<2, \quad t>0.  \]
The fundamental solution of this equation, denoted by $p_\alpha(t,x)$, can be computed via the Fourier transform by 
\[\widehat{p_\alpha}(t,\xi)= e^{-t\, \vert \xi \vert^\alpha}.\]
Moreover, in the spatial variable the fundamental solution $p_\alpha$ is given by 
\begin{equation}\label{pa}
  p_\alpha(t,x)= \frac{1}{t^{\frac{1}{\alpha}}}P_\alpha \left( \frac{x}{t^{\frac{1}{\alpha}}}\right),   
\end{equation}
where the function $P_\alpha$ is the inverse Fourier transform of $e^{-\vert \xi \vert^\alpha}$. See \cite[Chapter 13]{Jacob} for more details. It is well-known that for $1<\alpha<2$ the functions $P_\alpha$ are smooth and positive. In addition, they verify the following pointwise inequalities 
\[ 0<P_\alpha(x)\leq \frac{C}{(1+\vert x \vert)^{n+\alpha}}, \quad \vert \nabla P_\alpha(x)\vert \leq \frac{C}{(1+\vert x \vert^{n+\alpha+1})}, \] for a constant $C>0$ and for all $x\in \Rn$. These inequalities  allow us to derive the following estimates. 

\begin{Proposition}[$L^p-$estimates]\label{prop:1}
	For $1 \leq p \leq +\infty$, there exists a constant $C_{n,p} > 0$, which depends on the dimension $n \in \N^*$ and the parameter $p$, such that for every $1 < \alpha < 2$ and for every $t > 0$, we have
		\begin{enumerate}
			\item $\norm{p_{\alpha}(t, \cdot)}_{L^p} \leq C_{n,p} \ t^{-\frac{n}{\alpha}\left(1-\frac{1}{p}\right)}$,
			\item $\norm{\nabla p_{\alpha}(t, \cdot)}_{L^p} \leq C_{n,p} \ t^{-\frac{1+n\left(1-\nicefrac{1}{p}\right)}{\alpha}}$.
		\end{enumerate}
\end{Proposition}
Moreover we have:
\begin{Proposition}[$L^p-$continuity]\label{prop:2}
	Let $1 \leq p \leq +\infty$. For every $\varphi \in L^p(\R^n)$, we have
		\begin{align*}
			\lim_{t \to 0^{+}}\norm{p_{\alpha}(t, \cdot)*\varphi - \varphi}_{L^p}=0.
		\end{align*}
\end{Proposition}

On the other hand, by  the identity $\ds{\widehat{p_\alpha}(t,\xi)=e^{-t\, \vert \xi \vert^{\alpha}}}$  we have the following known results in the setting of the Sobolev spaces: 
\begin{Proposition}[$\dot{H}^s$ and $H^s$ estimates]\label{prop:3}
	Let $s_1, s_2 \geq 0$. Then, there is a constant $C_{n,s_2} > 0$, which depends on the dimension $n \in \N^*$ and the parameter $s_2$, such that for every $1 < \alpha \leq  2$ and for every $t > 0$, we have:
\begin{enumerate}
\item[1)]  $\ds{\norm{p_{\alpha}(t, \cdot)*\varphi}_{\dot{H}^{s_1 +s_2}} \leq C_{n,s_2} \ t^{-\frac{s_2}{\alpha}} \norm{\varphi}_{\dot{H}^{s_1}}.}$

\item[2)]  $\ds{\norm{p_\alpha(t, \cdot)*\varphi}_{H^{s_1 + s_2}} \leq C_{n,s_2} \left(1 + t^{-\nicefrac{s_2}{\alpha}}\right) \norm{\varphi}_{H^{s_1}}.}$
\end{enumerate}	
\end{Proposition}
\pv In order to verify the first point, we just write: 
\[ \Vert p_\alpha(t,\cdot)\ast \varphi \Vert^{2}_{\dot{H}^{s_1+s_2}}=\int_{\Rn}\vert \xi \vert^{2(s_1+s_2)} e^{-2 t \vert \xi \vert^\alpha}\vert\widehat{\varphi}(\xi)\vert^2 \, d\xi \leq t^{-\frac{2s_2}{\alpha}} \left( \sup_{\xi \in \Rn} \vert t^{1/\alpha} \xi \vert^{2 s_2}e^{-2 \vert t^{1/\alpha} \xi\vert^{\alpha}}\right) \int_{\Rn}\vert \xi \vert^{2 s_1}\vert \widehat{\varphi}(\xi)\vert^2\, d\xi.  \]

To verify the second point, let us start  by writing
\[ 	\norm{p_\alpha(t, \cdot) * \varphi}_{H^{s_1 + s_2}} = \norm{p_\alpha(t, \cdot) * \varphi}_{L^2} + \norm{p_\alpha(t, \cdot) * \varphi}_{\dot{H}^{s_1 + s_2}}. \] 
Then, for the first term on the right-hand side,  by the Young's inequalities and the point $1$  in Proposition \ref{prop:1}, we obtain 
\begin{align}\label{eq:9}
\norm{p_\alpha(t, \cdot) * \varphi}_{L^2} \leq \norm{p_\alpha(t, \cdot)}_{L^1} \norm{\varphi}_{L^2} \leq c \norm{\varphi}_{L^2} \leq c \norm{\varphi}_{H^{s_1}},
\end{align}
while for the second term on the right-hand side, by the point $1$ proven above we can write: 			
\begin{align}\label{eq:11}
\norm{p_\alpha(t, \cdot) * \varphi}_{\dot{H}^{s_1 + s_2}} \leq c_{n,s_2} t^{-\nicefrac{s_2}{\alpha}} \norm{\varphi}_{H^{s_1}}.
\end{align}
Thus, the desired estimate follows directly from \eqref{eq:9} and \eqref{eq:11}. \finpv

\begin{Proposition}[$H^s-$ and $\dot{H}^s-$continuity]\label{prop:4}
	Let $s_1, s_2 \geq 0$ and $\varepsilon > 0$. There exists a constant $C_{n,s_2, \varepsilon} > 0$, which depends on the dimension $n \in \N^*$, the parameters $s_2$ and $\varepsilon$, such that for every $1 < \alpha < 2$ and for every $t_1, t_2 > \varepsilon$, we have
		\begin{enumerate}
			\item $\norm{p_{\alpha}(t_1, \cdot)*\varphi - p_{\alpha}(t_2, \cdot)*\varphi}_{\dot{H}^{s_1 +s_2}} \leq C_{n,s_2,\varepsilon} \abs{t_1 - t_2}^{1/2} \norm{\varphi}_{\dot{H}^{s_1}}$,
			\item $\norm{p_{\alpha}(t_1, \cdot)*\varphi - p_{\alpha}(t_2, \cdot)*\varphi}_{H^{s_1 +s_2}} \leq C_{n,s_2,\varepsilon} \abs{t_1 - t_2}^{1/2} \norm{\varphi}_{H^{s_1}}$.\\
		\end{enumerate}
\end{Proposition}
\pv  To verify the first point, without loss of generality  we shall  assume that  $ t_1>t_2>\varepsilon$. Then we write 
\begin{equation*}
\begin{split}
&\Vert p_{\alpha}(t_1, \cdot)*\varphi - p_{\alpha}(t_2, \cdot)*\varphi \Vert^{2}_{\dot{H}^{s_1+s_2}}=\int_{\Rn} \vert \xi \vert^{2(s_1+s_2)} \vert e^{-t_1 \vert \xi \vert^\alpha} - e^{-t_2 \vert \xi \vert^\alpha}\vert^2 \vert \widehat{\varphi}(\xi)\vert^2\, d\xi  \\
=& \,  \int_{\Rn} \vert \xi \vert^{2 s_2} e^{-2 t_2 \vert \xi \vert^\alpha} \vert e^{-(t_1-t_2) \vert \xi \vert^\alpha} - 1\vert^2\, \vert \xi \vert^{2 s_1}  \vert\widehat{\varphi}(\xi)\vert^2\, d\xi\\
\leq&  \,  t^{-\frac{2s_2}{\alpha}}_{2} \,  \left( \sup_{\xi \in \Rn} \vert t^{1/\alpha}_{2} \xi \vert^{2 s_2}e^{- \vert t^{1/\alpha}_{2} \xi\vert^{\alpha}}\right)\, \int_{\Rn} e^{-t_2 \vert \xi \vert^\alpha} \vert e^{-(t_1-t_2) \vert \xi \vert^\alpha} - 1\vert^2\, \vert \xi \vert^{2 s_1}  \vert\widehat{\varphi}(\xi)\vert^2\, d\xi\\
\leq & \,  \varepsilon^{-\frac{2s_2}{\alpha}} C_{n,s_2} \, \int_{\Rn} e^{-\varepsilon \vert \xi \vert^\alpha} \vert e^{-(t_1-t_2) \vert \xi \vert^\alpha} - 1\vert^2\, \vert \xi \vert^{2 s_1}  \vert\widehat{\varphi}(\xi)\vert^2\, d\xi\\
\leq & \,  C_{n,s_2,\varepsilon}\, \int_{\Rn} e^{-\varepsilon \vert \xi \vert^\alpha} \vert e^{-(t_1-t_2) \vert \xi \vert^\alpha} - 1\vert^2\, \vert \xi \vert^{2 s_1}  \vert\widehat{\varphi}(\xi)\vert^2\, d\xi.
\end{split}    
\end{equation*}
We study now the expression $\ds{\vert e^{-(t_1-t_2) \vert \xi \vert^\alpha} - 1\vert^2}$. First, we remark that since we have $t_1>t_2$ then the expression $\ds{\vert e^{-(t_1-t_2) \vert \xi \vert^\alpha} - 1\vert}$ is uniformly  bounded and we can write 
\[ \vert e^{-(t_1-t_2) \vert \xi \vert^\alpha} - 1\vert^2=\vert e^{-(t_1-t_2) \vert \xi \vert^\alpha} - 1\vert\, \vert e^{-(t_1-t_2) \vert \xi \vert^\alpha} - 1\vert\leq C \vert e^{-(t_1-t_2) \vert \xi \vert^\alpha} - 1\vert.  \]
Now, by the mean value theorem in the temporal variable we have $\ds{\vert e^{-(t_1-t_2) \vert \xi \vert^\alpha} - 1\vert \leq C\, \vert \xi \vert^\alpha \, \vert t_1-t_2 \vert}$. Thus, gathering these estimates we  get 
\[ \vert e^{-(t_1-t_2) \vert \xi \vert^\alpha} - 1\vert^2 \leq C\, \vert \xi \vert^{\alpha}\vert t_1-t_2\vert.\]
Getting back to the last integral we finally have:
\begin{equation*}
\begin{split}
 \Vert p_{\alpha}(t_1, \cdot)*\varphi - p_{\alpha}(t_2, \cdot)*\varphi \Vert^{2}_{\dot{H}^{s_1+s_2}} \leq &  C_{n,s_2,\varepsilon}\,\vert t_1-t_2 \vert \, \int_{\Rn} e^{-\varepsilon \vert \xi \vert^\alpha}\vert \xi \vert^\alpha \, \vert \xi \vert^{2 s_1}  \vert\widehat{\varphi}(\xi)\vert^2\, d\xi \\
 \leq & C_{n,s_2,\varepsilon}\,\vert t_1-t_2 \vert \left(\sup_{\xi \in \Rn}  e^{-\varepsilon \vert \xi \vert^\alpha}\vert \xi \vert^{\alpha}\right) \Vert \varphi \Vert^{2}_{\dot{H}^{s_1}}\\
 \leq &  C_{n,s_2,\varepsilon}\,\vert t_1-t_2 \vert\, \Vert \varphi \Vert^{2}_{\dot{H}^{s_1}},
\end{split}    
\end{equation*}
hence, the first point is verified. The second point essentially follows these sames lines.\finpv

%%%%%%%%%%%%%%%%%%%%%%%%%%%%%%%%%%%%%%%%%%%%%%%%%%%%%

\section{Global well-posedness and regularity: proof of  Theorem \ref{Th-WP}}\label{Sec:WP-Reg}
Let $1 < \alpha \leq 2$ fixed, and let $u_{0} \in L^1\cap H^s(\Rn)$ be an initial datum. The result stated in Theorem \ref{Th-WP} is well-known for the case $\alpha=2$, see for instance \cite{Biler1}, \cite{Dro} and \cite{Dro2}. Consequently, we just consider the range $1<\alpha <2$.  As mentioned, the proof of this theorem is rather standard  but, for the reader's convenience, we shall detail some technical estimates. 

\medskip 

%%%%%%%% STEP 1 %%%%%%%%%%%
\textbf{Step 1: Local well-posedness}.  We consider the (equivalent) mild formulation given in (\ref{u-alpha}), where  the nonlinear term defines a multi-linear form in the variable $u$ (see the expression (\ref{M(u)}) below). In order to construct a solution of the equation (\ref{u-alpha}) we will use the  \emph{Picard's contraction principle} for a time $0<T<+\infty$ small enough. We thus consider  the Banach space 
\begin{align}\label{Espace}
E_T = \mathcal{C}\big([0, T], L^1(\R^n)\big) \cap \mathcal{C}\big([0, T], H^{s}(\R^n)\big),
\end{align}
endowed with the norm
\begin{align}\label{Norm}
\norm{u}_{E_T} = \sup_{0 \leq t \leq T} \norm{u(t, \cdot)}_{L^1} + \sup_{0 \leq t \leq T} \norm{u(t, \cdot)}_{H^{s}}.
\end{align}
Then, we will prove the following::
\begin{Proposition}\label{Th-tech-1} Let  $s>n/2$ and let $u_{0} \in L^1\cap H^s(\Rn)$ be an initial datum.  Moreover, let $1<\alpha<2$. Then, there exists a time  given by: 
	\begin{align}\label{Time-T}
		T = \dfrac{1}{2} \left[\dfrac{1-\frac{1}{\alpha}}{2^b c \abs{\eta} \Big(\norm{u_{0}}_{L^1} + \norm{u_{0}}_{H^s} \Big)^{b-1}}\right]^{\frac{\alpha}{\alpha-1}},
	\end{align}
where $c>0$ is a  numerical constant, and  oreover, there exists a function $u_\alpha \in E_T$ which is a solution of the equation (\ref{u-alpha}).  
\end{Proposition}
\pv We start by estimating the linear term in  the equation (\ref{u-alpha}).
\begin{Lemme}\label{Prop-Lin} Let $p_\alpha(t,x)$ be the kernel given in (\ref{pa}). Then   we have  $\Vert p_\alpha(t,\cdot)\ast u_{0,\alpha} \Vert_{E_T} \leq c\, (\Vert u_{0} \Vert_{L^1}+\Vert u_{0} \Vert_{H^s})$. \end{Lemme} 
\pv   We first observe that, due to Proposition \ref{prop:2} and the first point in Proposition \ref{prop:3}, the quantities $\Vert p_\alpha(t,\cdot)\ast u_{0} \Vert_{L^1}$ and $\Vert p_\alpha(t,\cdot)\ast u_{0} \Vert_{H^s}$ are continuous in the temporal variable.

\medskip 

On the other hand, by the  Young's inequalities and the point $1$ in Proposition \ref{prop:1} (with $p=1$)  we write
\[ \Vert p_\alpha(t,\cdot)\ast u_{0} \Vert_{L^1} \leq \norm{p_{\alpha}(t, \cdot)}_{L^1} \norm{u_{0}}_{L^1} \leq c\norm{u_{0}}_{L^1}. \]
We also write 
\[ \Vert p_\alpha(t,\cdot)\ast u_{0}\Vert_{H^s}\leq \Vert \widehat{p_\alpha}(t,\cdot)\Vert_{L^{\infty}}\Vert u_{0} \Vert_{H^s} \leq \Vert p_\alpha(t,\cdot)\Vert_{L^1}\Vert u_{0} \Vert_{H^s}\leq c \Vert u_{0} \Vert_{H^s}.
\]
to obtain the wished estimate. \finpv 

\medskip 

We study now  the nonlinear term in the equation (\ref{u-alpha}). For $b\in \mathbb{N}$ with $b\geq  2$,  we denote the multi-linear form  
\begin{equation}\label{M(u)}
M_b(u)= \int_{0}^{t}p_\alpha(t-\tau,\cdot)\ast \eta \cdot \nabla(u^{b})(\tau,\cdot)d\tau, 
\end{equation}
where, to simplify our writing, we have written the function $u$ instead of $u_\alpha$.  Then, we have the following  estimate 
\begin{Lemme}\label{Prop-Non-Lin} For $u\in E_T$ we have $M_b(u) \in E_T$. Moreover, the following estimate holds:\[ \ds{\Vert M_b(u)\Vert_{E_T} \leq c\,\vert \eta \vert\, \frac{T^{1-\nicefrac{1}{\alpha}}}{1-\nicefrac{1}{\alpha}}\, \Vert u \Vert^{b}_{E_T}}.\]  
\end{Lemme}
\pv   By  \textbf{\cite{Biler1}} we have $M_b(u)\in \mathcal{C}\big([0,T], L^1(\R^n)\big)$,  so it  remains to prove that  $M_b(u) \in \mathcal{C}\big([0,T], H^s(\R^n)\big)$. Indeed,  let $t_1, t_2 > 0$ and without loss of generality we assume that $0 < t_1 < t_2 \leq T$. Then  we write 
\begin{equation}\label{eqWP:7}
\begin{split}
 &\norm{\int_{0}^{t_1} p_{\alpha}(t_1 - \tau, \cdot) * \eta \cdot \nabla (u^b)(\tau, \cdot)  d\tau - \int_{0}^{t_2} p_{\alpha}(t_2 - \tau, \cdot) \ast \eta \cdot \nabla (u^b)(\tau, \cdot)  d\tau}_{H^s}\\
\leq & \left\Vert \int_{0}^{t_1} p_{\alpha}(t_1 - \tau, \cdot) \ast \eta \cdot \nabla (u^b)(\tau, \cdot)  d\tau - \int_{0}^{t_1} p_{\alpha}(t_2 - \tau, \cdot) \ast \eta \cdot \nabla (u^b)(\tau, \cdot)  d\tau \right\Vert_{H^s} \\
&+ \left\Vert \int_{0}^{t_1} p_{\alpha}(t_2 - \tau, \cdot) \ast \eta \cdot \nabla (u^b)(\tau, \cdot)  d\tau - \int_{0}^{t_2} p_{\alpha}(t_2 - \tau, \cdot) \ast \eta \cdot \nabla (u^b)(\tau, \cdot)  d\tau \right\Vert_{H^s}\\ 
\leq & \int_{0}^{t_1} \norm{p_{\alpha}(t_1 - \tau, \cdot) \ast \eta  u^b(\tau, \cdot) - p_{\alpha}(t_2 - \tau, \cdot) \ast \eta u^b(\tau, \cdot)}_{H^{s+1}} d\tau \\
&+ \int_{t_1}^{t_2} \norm{\nabla p_{\alpha}(t_2 - \tau, \cdot) \ast \eta  u^b(\tau, \cdot)}_{H^s} d\tau \\
= & R_{\alpha,1}(t_1,t_2) + R_{\alpha,2}(t_1,t_2). 
\end{split}    
\end{equation} 

For the first term on the right-hand side, by the point $2$ in Proposition \ref{prop:4} (with $s_1=s$ and $s_2=1$)  and  as $s>n/2$,  by the product laws in the Sobolev spaces we obtain: 
\begin{equation*}
\begin{split}
R_{\alpha,1}(t_1,t_2) \leq  &\,  c\, \int_{0}^{t_1}  \abs{t_1 - t_2}^{1/2} \abs{\eta} \norm{u^b(\tau, \cdot)}_{H^{s}} \ d\tau\leq c \abs{\eta} \abs{t_1 - t_2}^{1/2} \int_{0}^{t_1} \Vert u(\tau,\cdot)\Vert^{b}_{H^s} d \tau  \\
\leq  &\,  c \abs{\eta} \abs{t_1 - t_2}^{1/2} \, T \Vert u \Vert^{b}_{E_T}. 
\end{split}
\end{equation*}

Hence, we have $\ds{\lim_{t_1 \to t_2} R_{\alpha,1}(t_1,t_2)=0}$.  On the other hand, for the second term on the right-hand  side  we write 
\begin{equation*}
\begin{split}
R_{\alpha,2}(t_1,t_2)&= \int_{t_1}^{t_2} \norm{\nabla p_{\alpha}(t_2 - \tau, \cdot) \ast \eta \, u^b(\tau, \cdot)}_{L^2} + \norm{\nabla p_{\alpha}(t_2 - \tau, \cdot) \ast \eta \ u^b(\tau, \cdot)}_{\dot{H}^s} \ d\tau\\
&\leq \vert \eta \vert  \int_{t_1}^{t_2} \norm{\nabla p_{\alpha}(t_2 - \tau, \cdot)}_{L^1} \norm{ u^b(\tau, \cdot)}_{L^2}\, d\tau  + \vert \eta \vert \int_{t_1}^{t_1} \norm{p_{\alpha}(t_2 - \tau, \cdot) \ast   u^b(\tau, \cdot)}_{\dot{H}^{s+1}} \, d\tau\\
&= R_{\alpha,2,1}(t_1, t_2)+R_{\alpha,2,2}(t_1, t_2). 
\end{split}    
\end{equation*}
In order to estimate the term $R_{\alpha,2,1}(t_1,t_2)$, by the  H\"older inequalities, the second point in Proposition \ref{prop:1}, and moreover, the product laws in the Sobolev spaces, we write: 
\begin{equation*}
\begin{split}
 R_{\alpha,2,1}(t_1, t_2) \leq & c\, \vert \eta \vert \int_{t_1}^{t_2} (t_2 -\tau)^{-1/\alpha}\Vert u^b(\tau,\cdot)\Vert_{L^2}\, d \tau \leq   c\, \vert \eta \vert \int_{t_1}^{t_2} (t_2 -\tau)^{-1/\alpha}\Vert u^b(\tau,\cdot)\Vert_{H^s}\, d \tau\\
 \leq & c \abs{\eta} \norm{u}_{E_T}^b \dfrac{\abs{t_2 - t_1}^{1-\nicefrac{1}{\alpha}}}{1-\nicefrac{1}{\alpha}}. \end{split}    
\end{equation*}
In addition, in order to estimate $R_{\alpha,2,2}(t_1,t_2)$, by Proposition \ref{prop:3} (with $s_1=s$ and $s_2=1$), and  by using again the product laws in Sobolev spaces, we can write 
\begin{equation*}
R_{\alpha,2,2}(t_1, t_2) \leq c\, \vert \eta \vert \int_{t_1}^{t_2} (t_2-\tau)^{-1/\alpha} \Vert u^b(\tau,\cdot)\Vert_{\dot{H}^s} d\tau  \leq c \abs{\eta} \norm{u}_{E_T}^b \dfrac{\abs{t_2 - t_1}^{1-\nicefrac{1}{\alpha}}}{1-\nicefrac{1}{\alpha}}. \end{equation*}
By gathering the estimates made for the terms $R_{\alpha,2,1}(t_1,t_2)$ and $R_{\alpha,2,2}(t_1,t_2)$, we obtain $\ds{\lim_{t_1 \to t_2} R_{\alpha,2}(t_1,t_2)=0}$. We thus have $M_b(u)\in \mathcal{C}((0,T],H^s(\Rn))$. Now, we must prove the continuity at $t=0$. For this we will verify the estimate

\begin{equation}\label{estim}
\left\Vert \int_{0}^{t} p_\alpha(t-\tau, \cdot)\ast \eta \cdot \nabla (u^b)(\tau, \cdot)d \tau \right\Vert_{H^s} \leq c\vert \eta \vert \Vert u \Vert^{b}_{E_T}\frac{t^{1-1/\alpha}}{1-1/ \alpha}.     
\end{equation}

Indeed, by the Young's inequalities, the second point in Proposition \ref{prop:1}, Proposition \ref{prop:3}, and moreover, the product laws in Sobolev spaces we can write: 
\begin{align*}
	&\norm{\int_{0}^{t} p_{\alpha}(t - \tau, \cdot) * \eta \cdot \nabla (u^b)(\tau, \cdot) \ d\tau}_{H^s}\\
	%&\leq \int_{0}^{t} \norm{\nabla p_{\alpha}(t - \tau, \cdot) * \eta \ u_{\alpha}^b(\tau, \cdot)}_{H^s} \ d\tau\\
	%&\leq \int_{0}^{t} \norm{\nabla p_{\alpha}(t - \tau, \cdot) * \eta \ u_{\alpha}^b(\tau, \cdot)}_{L^2} + \norm{\nabla p_{\alpha}(t - \tau, \cdot) * \eta \ u_{\alpha}^b(\tau, \cdot)}_{\dot{H}^s} \ d\tau\\
	\leq &\int_{0}^{t} \norm{\nabla p_{\alpha}(t - \tau, \cdot)}_{L^1} \norm{\eta \ u^b(\tau, \cdot)}_{L^2} + \norm{p_{\alpha}(t - \tau, \cdot) * \eta \ u^b(\tau, \cdot)}_{\dot{H}^{s+1}}\ d\tau\\
	%\end{align*}
	%Thus, by $L^p$ and $\dot{H}^s-$estimates of $p_\alpha(t, \cdot)$:
	%\begin{align*}
	%&\norm{\int_{0}^{t} p_{\alpha}(t - \tau, \cdot) * \eta \cdot \nabla (u_{\alpha}^b)(\tau, \cdot) \ d\tau}_{H^s}\\
	\leq &\int_{0}^{t} c (t - \tau)^{-\nicefrac{1}{\alpha}} \abs{\eta} \norm{u^b(\tau, \cdot)}_{L^2} + c (t - \tau)^{-\nicefrac{1}{\alpha}} \abs{\eta} \norm{u^b(\tau, \cdot)}_{\dot{H}^s}\ d\tau\\
	%& \leq c_n \abs{\eta} \int_{0}^{t} (t - \tau)^{-\nicefrac{1}{\alpha}} \norm{u_{\alpha}^b (\tau, \cdot)}_{H^s}\ d\tau\\
	\leq &\ c \abs{\eta} \norm{u}_{E_T}^b \int_{0}^{t} (t - \tau)^{-\nicefrac{1}{\alpha}} \ d\tau \leq  c\vert \eta \vert \Vert u \Vert^{b}_{E_T}\frac{t^{1-1/\alpha}}{1-1/ \alpha}. %\\
	%& \leq c_n \abs{\eta} \norm{u_\alpha}_{E_T}^b \dfrac{t^{1-\nicefrac{1}{\alpha}}}{1-\nicefrac{1}{\alpha}}.
	\end{align*}
	
Once we have $M_b(u) \in E_T$, 	we  verify now the estimate stated in Lemma \ref{Prop-Non-Lin}. First  note that by the estimate (\ref{estim}) we can  write 
\begin{equation}\label{estim1}
  \sup_{t \in [0,T]}\, \left\Vert \int_{0}^{t} p_\alpha(t-\tau, \cdot)\ast \eta \cdot \nabla (u^b)(\tau, \cdot)d \tau \right\Vert_{H^s} \leq c\vert \eta \vert \frac{T^{1-1/\alpha}}{1-1/ \alpha}\, \Vert u \Vert^{b}_{E_T}.    
\end{equation}
On the other hand, by  applying  the Young inequalities and the point $2$ in Proposition \ref{prop:1} we have 
\begin{equation*}
\begin{split}
& \norm{\int_{0}^{t} p_{\alpha}(t - \tau, \cdot) \ast  \eta \cdot \nabla (u^b)(\tau, \cdot)  d\tau}_{L^1}
	\leq  \int_{0}^{t} \norm{\nabla p_{\alpha}(t - \tau, \cdot) \ast \eta \ u^b(\tau, \cdot)}_{L^1}  d\tau\\
	%& \leq \int_{0}^{t} \norm{\nabla p_{\alpha}(t - \tau, \cdot)}_{L^1} \norm{\eta\ u_{\alpha}^b(\tau, \cdot)}_{L^1} \ d\tau\\
	 \leq & c \abs{\eta}\, \int_{0}^{t} (t - \tau)^{-\nicefrac{1}{\alpha}}  \norm{u^b(\tau, \cdot)}_{L^1}  d\tau  \leq   c \abs{\eta}\, \int_{0}^{t} (t - \tau)^{-\nicefrac{1}{\alpha}}  \Vert u(\tau,\cdot)\Vert^{b-1}_{L^{\infty}} \norm{u(\tau, \cdot)}_{L^1}  d\tau.  
\end{split}    
\end{equation*}
Since $s>n/2$ we have the embedding $H^s(\Rn)\subset L^{\infty}(\Rn)$, and thus we can write 
\begin{equation*}
\begin{split}
 &c \abs{\eta}\, \int_{0}^{t} (t - \tau)^{-\nicefrac{1}{\alpha}}  \Vert u(\tau,\cdot)\Vert^{b-1}_{L^{\infty}} \norm{u(\tau, \cdot)}_{L^1}  d\tau \leq    c \abs{\eta}\, \int_{0}^{t} (t - \tau)^{-\nicefrac{1}{\alpha}}  \Vert u(\tau,\cdot)\Vert^{b-1}_{H^s} \norm{u(\tau, \cdot)}_{L^1}  d\tau\\
 \leq & c\vert \eta \vert \left( \sup_{\tau \in [0,T]} \Vert u(\tau,\cdot)\Vert^{b-1}_{H^s} \right) \left( \sup_{\tau \in [0,T]} \Vert u(\tau,\cdot)\Vert_{L^1} \right) \frac{T^{1-1/\alpha}}{1-1/ \alpha} \leq c\vert \eta \vert  \,\frac{T^{1-1/\alpha}}{1-1/ \alpha}\, \Vert u \Vert^{b}_{E_T}.
\end{split}    
\end{equation*}
Then, we have 
\begin{equation}\label{estim2}
\begin{split}
 \sup_{t \in [0,T]}   \norm{\int_{0}^{t} p_{\alpha}(t - \tau, \cdot) \ast  \eta \cdot \nabla (u^b)(\tau, \cdot)  d\tau}_{L^1}
	\leq  &\,  \int_{0}^{t} \norm{\nabla p_{\alpha}(t - \tau, \cdot) \ast \eta \ u^b(\tau, \cdot)}_{L^1}  d\tau \\
	 \leq  &\,   c\vert \eta \vert  \,\frac{T^{1-1/\alpha}}{1-1/ \alpha}\, \Vert u \Vert^{b}_{E_T}.
\end{split}	
\end{equation}
Finally, by (\ref{estim1}) and (\ref{estim2}) we obtain the desired estimate. This lemma is proven. \finpv

\medskip 

Once we have  Lemmas \ref{Prop-Lin} and \ref{Prop-Non-Lin} at our disposal,  the rest of the proof of Proposition  \ref{Th-tech-1} follows from standard arguments. \finpv. 

\medskip

\textbf{Step 2: Regularity}. We define the space $H^{\infty}(\Rn)$ as $\ds{H^{\infty}(\R^n) = \bigcap_{s \geq 0} H^{s}(\R^n)}$. 
\begin{Proposition}\label{prop:6}
	Let $u_\alpha \in E_{T}$ be the unique solution of the equation (\ref{u-alpha}) given by Theorem \ref{Th-tech-1}. This solution satisfies $u_\alpha \in \mathcal{C}\big((0,T],H^{\infty}(\R^n)\big)$.  Moreover, we have $u_\alpha \in \mathcal{C}^{1}((0,T], \mathcal{C}^{\infty}(\Rn))$; and for $0<t\leq T$  the solution $u_\alpha$  verifies the differential equation (\ref{Eq3}) in the classical sense.    
	\end{Proposition}  

\pv We will verify that each term on the right-hand side  in the equation  \eqref{u-alpha} belongs to the space $\mathcal{C}\big([0,T],H^{\infty}(\R^n)\big)$. For the first (linear) term,  by the second point in Proposition  \ref{prop:3}, and moreover, by the second point in  Proposition \ref{prop:4}, we directly have  $\ds{p_\alpha*u_{0,\alpha}\in \mathcal{C}\big((0,T],H^{\infty}(\R^n)\big)}$.

\medskip

For the second (nonlinear) term,  we recall that  by (\ref{estim})  for all time $0<t\leq T$  we have $\int_{0}^{t} p_\alpha(t-\tau, \cdot)\ast \eta \cdot \nabla (u^{b}_{\alpha})(\tau, \cdot)d \tau \in H^s(\R^n)$.  Then,  we will prove that for  $\sigma > 0$ small enough  we also have: $\int_{0}^{t} p_\alpha(t-\tau, \cdot)\ast \eta \cdot \nabla (u^{b}_{\alpha})(\tau, \cdot)d \tau \in H^{s+\sigma}(\R^n)$. Indeed, by using the second point in Proposition  \ref{prop:3}, for $\sigma >0$ (which we shall  set later)  we write 
\begin{equation*}
\begin{split}
&\, \norm{\int_{0}^{t} p_\alpha(t-\tau, \cdot)\ast \eta \cdot \nabla (u^{b}_{\alpha})(\tau, \cdot)d \tau}_{H^{s+\sigma}}  \leq \, c\, \vert \eta\vert \,  \int_{0}^{t} \norm{p_\alpha(t-\tau,\cdot) \ast  u^{b}_{\alpha} (\tau,\cdot)}_{H^{s +\sigma + 1}}\ d\tau\\
&\leq c\, \vert \eta\vert \, \int_{0}^{t} \left[1 + (t-\tau)^{-\nicefrac{(\sigma + 1)}{\alpha}}\right] \norm{ u_\alpha^b (\tau,\cdot)}_{H^s}\ d\tau \leq c\, \vert \eta\vert \,  \norm{u_\alpha}^b_{E_{T}} \int_{0}^{t}\left[  1 + (t-\tau)^{-\nicefrac{(\sigma + 1)}{\alpha}}\right] d\tau.
\end{split}
\end{equation*}
We thus set $0<\sigma < \alpha -1$ (recall that we have $1<\alpha < 2$) to obtain that the last integral above computes down as 
\[  \int_{0}^{t} 1 + (t-\tau)^{-\nicefrac{(\sigma + 1)}{\alpha}}\ d\tau = t + \dfrac{t^{1- \nicefrac{(\sigma + 1)}{\alpha}}}{1- \nicefrac{(\sigma + 1)}{\alpha}}.\]
Then, for all time $0<t\leq T$ we  obtain the estimate: 
\[ \norm{\int_{0}^{t} p_\alpha(t-\tau, \cdot)\ast \eta \cdot \nabla (u^{b}_{\alpha})(\tau, \cdot)d \tau}_{H^{s+\sigma}} \leq c \abs{\eta} \norm{u_\alpha}^b_{E_{T}} \left[t + \dfrac{t^{1- \nicefrac{(\sigma + 1)}{\alpha}}}{1- \nicefrac{(\sigma + 1)}{\alpha}}\right]. \]

We will show now that we have $\ds{\int_{0}^{t} p_\alpha(t-\tau, \cdot)\ast \eta \cdot \nabla (u^{b}_{\alpha})(\tau, \cdot)d \tau \in \mathcal{C}\big((0,T], H^{s+\sigma}(\R^n)\big)}$.  Let $0 < t_1, t_2 < T$, where, always without loss of generality  we shall assume  that $t_1 < t_2$.  Then we write: 
\begin{align}\label{eq:12}
		&\norm{\int_{0}^{t_2}  p_\alpha(t_2-\tau,\cdot) \ast \eta \cdot \nabla (u_\alpha^b) (\tau,\cdot)\ d\tau - \int_{0}^{t_1}  p_\alpha(t_1-\tau,\cdot) \ast  \eta \cdot \nabla  (u_\alpha^b)(\tau,\cdot)\ d\tau}_{H^{s+\sigma}}\notag\\
		\leq &\norm{\int_{0}^{t_2}  p_\alpha(t_2-\tau,\cdot) \ast  \eta\cdot \nabla  (u_\alpha^b)(\tau,\cdot)\ d\tau - \int_{0}^{t_1}  p_\alpha(t_2-\tau,\cdot) \ast \eta \cdot \nabla (u_\alpha^b)(\tau,\cdot)\ d\tau}_{H^{s+\sigma}}\notag\\
		&+ \norm{\int_{0}^{t_1}  p_\alpha(t_2-\tau,\cdot) \ast \eta \cdot \nabla (u_\alpha^b)(\tau,\cdot)\ d\tau - \int_{0}^{t_1}  p_\alpha(t_1-\tau,\cdot) \ast  \eta \cdot \nabla (u_\alpha^b)(\tau,\cdot)\ d\tau}_{H^{s+\sigma}}\notag\\
		= & \ \tilde{R}_{\alpha,1}(t_1,t_2) + \tilde{R}_{\alpha,2}(t_1,t_2),
	\end{align}
where, we must study the terms $\tilde{R}_{\alpha,1}(t_1,t_2)$ and $\tilde{R}_{\alpha,2}(t_1,t_2)$. For the term $\tilde{R}_{\alpha,1}(t_1,t_2)$, by the second point in Proposition \ref{prop:3} we can write: 
\begin{align*}
\tilde{R}_{\alpha,1}(t_1,t_2) \leq \int_{t_1}^{t_2} \norm{p_\alpha(t_2-\tau,\cdot) * \eta\ (u_\alpha^b)(\tau,\cdot)}_{H^{s+\sigma + 1}} d\tau \leq C \int_{t_1}^{t_2} \left[1 + (t_2 - \tau)^{-\nicefrac{(\sigma + 1)}{\alpha}}\right] \norm{\eta\ (u_\alpha^b)(\tau,\cdot)}_{H^s} d\tau.
\end{align*}
Since $0<\sigma < \alpha -1$  the integral above computes down as 
\[ \displaystyle \int_{t_1}^{t_2} 1 + (t_2-\tau)^{-\nicefrac{(\sigma + 1)}{\alpha}}\ d\tau = (t_2 - t_1) + \dfrac{(t_2 -t_1)^{1-{\nicefrac{(\sigma + 1)}{\alpha}}}}{{1-{\nicefrac{(\sigma + 1)}{\alpha}}}}. \]
Hence, we have: 
\begin{align}\label{eq:13}
\tilde{R}_{\alpha,1}(t_1,t_2) \leq \ c \abs{\eta} \norm{u}^{b}_{E_{T}} \left[(t_2 - t_1) + \dfrac{(t_2 -t_1)^{1-{\nicefrac{(\sigma + 1)}{\alpha}}}}{{1-{\nicefrac{(\sigma + 1)}{\alpha}}}}\right].
\end{align} 
For the term  $\tilde{R}_{\alpha,2}(t_1,t_2)$, always by the second point in Proposition \ref{prop:4}, we can write:
 \begin{align}\label{eq:14}
\tilde{R}_{\alpha,2}(t_1,t_2) %& \norm{\int_{0}^{t_1} \nabla p_\alpha(t_2-\tau,\cdot) * \eta\ (u_\alpha^b)(\tau,\cdot)\ d\tau - \int_{0}^{t_1} \nabla p_\alpha(t_1-\tau,\cdot) * \eta\ (u_\alpha^b)(\tau,\cdot)\ d\tau}_{H^{s+\sigma}} \\
%& \leq \int_{0}^{t_1} \norm{\nabla \left[p_\alpha(t_2-\tau, \cdot) * \eta\ (u_\alpha^b)(\tau,\cdot)- p_\alpha(t_1-\tau, \cdot) * \eta\ (u_\alpha^b)(\tau,\cdot) \right]}_{H^{s+\sigma}}\ d \tau  \\
\leq &c\,\vert \eta \vert, \int_{0}^{t_1}  \norm{p_\alpha(t_2-\tau, \cdot)  \eta\ (u_\alpha^b)(\tau,\cdot) - p_\alpha(t_1-\tau, \cdot) \ast   (u_\alpha^b)(\tau,\cdot)}_{H^{s+\sigma + 1}}\ d \tau\notag\\
\leq &\ c\, \vert \eta \vert\,  \abs{t_1-t_2}^{1/2}\, \int_{0}^{t_1} \norm{ u_\alpha^b(\tau,\cdot)}_{H^s}\ d \tau \leq c\, \vert \eta \vert\,  \abs{t_1-t_2}^{1/2}\,  T\, \norm{u}^{b}_{E_{T}}.
\end{align}

Therefore, for $0 < \sigma < \alpha-1$, by \eqref{eq:13} and \eqref{eq:14} we have
	\begin{align*}
		\int_{0}^{t}  p_\alpha(t-\tau,\cdot) \ast \eta \cdot \nabla  (u_\alpha^b)(\tau,\cdot)\ d\tau \in \mathcal{C}\big((0,T], H^{s+\sigma}(\R^n)\big).
	\end{align*}

 At this point, we have proven that $\ds{u_\alpha \in  \mathcal{C}((0, T_0], H^{s+\delta}(\Rn))}$ and by repeating this process (in order to obtain a gain of regularity for the nonlinear term) we conclude that $\ds{u_\alpha \in  \mathcal{C}((0, T], H^{\infty}(\Rn))}$.  
 
\medskip

With this information at our disposal,  we can verify now that  for all $0<t\leq T$ and for all  multi-index ${\bf a} \in \mathbb{N}^n$   we have  $\partial^{\bf a}_{x}u_\alpha(t,\cdot) \in \mathcal{C}((0,T], \mathcal{C} \cap L^{\infty}(\Rn))$. Indeed, let ${\bf a}=(a_1,\cdots, a_n) \in \mathbb{N}^n$ be a multi-index, where we denote by $\vert {\bf a} \vert = a_1 + \cdots + a_n$ its size. Then, for  $\frac{n}{2} < s_1 < \frac{n}{2}+1$ we set $s= \vert {\bf a} \vert + s_1$. Since we have $\ds{u_\alpha \in  \mathcal{C}((0, T], H^{\infty}(\Rn))}$ then we get $\ds{\partial^{{\bf a}}_{x} u_\alpha(t,\cdot) \in H^{s_1}}(\Rn)$. Moreover, since  $\frac{n}{2}<s_1$ we have the continuous embedding $H^{s_1}(\Rn)\subset L^{\infty}(\Rn)$, hence we conclude that $\partial^{{\bf a}}_{x} u_\alpha(t,\cdot) \in L^{\infty}(\Rn)$. 

\medskip 

On the other hand, we recall that we have the identification $H^{s_1}(\Rn)= B^{s_1}_{2,2}(\Rn)$ (where $B^{s_1}_{2,2}(\Rn)$ denotes a non-homogeneous  Besov space \cite{BahouriDanchinCheman}). Moreover, we also have the continuous embedding $B^{s_1}_{2,2}(\Rn)\subset B^{s_1 - n/1}_{\infty, \infty}(\Rn)\subset \dot{B}^{s_1-n/2}_{\infty,\infty}(\Rn)$.

\medskip

We thus have  $\ds{\partial^{{\bf a}}_{x} u_\alpha(t,\cdot) \in \dot{B}^{s_1-n/2}_{\infty,\infty}(\Rn)}$. But, since $\frac{n}{2} < s_1 < \frac{n}{2}+1$ then we have $0< s_1 -\frac{n}{2}<1$, and thereafter, by definition of the homogeneous Besov space  $\dot{B}^{s_1-n/2}_{\infty,\infty}(\Rn)$ (see always \cite{BahouriDanchinCheman}) we get that  $\partial^{{\bf a}}_{x} u_\alpha(t,\cdot)$ is a $\beta-$H\"older continuous functions with parameter $\beta= s_1 - \frac{n}{2}\in (0,1)$. 

\medskip

We  have  proven that $u_\alpha \in \mathcal{C}((0,T], \mathcal{C}^{\infty}(\Rn))$. and we write $\partial_t u_\alpha = - (-\Delta)^{\alpha / 2} u_\alpha - \eta \cdot \nabla (u^{b}_{\alpha})$ to obtain that $\partial_t u_\alpha \in \mathcal{C}((0,T], \mathcal{C}^{\infty}(\Rn))$. Finally, we conclude that $u_\alpha \in \mathcal{C}^{1}((0, T], \mathcal{C}^{\infty}(\Rn))$.  Proposition \ref{prop:6} is proven.  \finpv 
%%%%%%%%%%%%%%%%%%%%%%%%%%%%%%%%%%%%%%%%%%%%%%

\medskip

\textbf{Step 4: Global in time existence}. By following similar arguments of \cite{Cui} (see the proof of Theorem 2, page 9) we have the following result.
\begin{Proposition}\label{Th-Tech-4} Let $u_0 \in L^1\cap H^s(\R^n)$ be an initial data  and let  $T^* > 0$ be the maximal time of existence of the unique corresponding arising  solution $u_\alpha \in E_{T^*}$ (given by Theorem \ref{Th-tech-1}) to the problem \eqref{u-alpha}. Then we have $T^* = +\infty$. 
\end{Proposition} 
\pv  Let us briefly explain the general idea of the proof.  We assume that  $T^* < +\infty$. Then we will extend the solution  $u_\alpha$  beyond the time $T^{*}$ to obtain a contradiction.  We thus  conclude that $T^* = +\infty$. 
 
\medskip

We start by defining the  following function  $\func{T(\cdot)}{[0, +\infty[}{[0, +\infty[}$   such that for each initial datum $w_{0} \in L^1 \cap H^s(\R^n)$ the quantity $T\left(\norm{w_{0}}_{L^1}\right)$ is given by the expression  
	\begin{align*}
		T\left(\norm{w_{0}}_{L^1}\right) = \dfrac{1}{2} \left[\dfrac{1-\nicefrac{1}{\alpha}}{2^b c \abs{\eta} \left(\norm{w_{0}}_{L^1} + \norm{w_{0}}_{H^s} \right)^{b-1}}\right]^{\frac{\alpha}{\alpha-1}}.
	\end{align*}
We  recall that $T\left(\norm{w_{0}}_{L^1}\right)$ is precisely the  first time of the existence of the solution $w_\alpha$  to the equation \eqref{u-alpha}, which is given by Theorem \ref{Th-tech-1}.   Additionally, the \emph{key remark} is that  this function is decreasing in the variable $\norm{w_{0}}_{L^1}$.

\medskip  

On the other hand, by \textbf{\cite{Biler1}} we known that  for every initial datum $w_{0} \in L^1 \cap H^s(\R^n)$ we have a  unique solution $w_\alpha \in \mathcal{C}\big([0,+\infty[,L^1(\R^n)\big)$ of the equation \eqref{u-alpha}. Moreover, for every time $t > 0$ we have the estimate
	\begin{align}\label{eq:15}
		\norm{w_\alpha(t, \cdot)}_{L^1} \leq \norm{w_{0}}_{L^1}.
	\end{align}
Since the function $T(\cdot)$  defined above is decreasing in the variable $\norm{w_{0}}_{L^1}$,  for the initial datum $u_{0} \in L^1\cap H^s(\R^n)$ we can set a  time $0<T_1 < T^{*}$ such that for all $w_{0} \in L^1\cap H^s(\R^n)$ with $\norm{w_{0}}_{L^1} \leq \norm{u_{0}}_{L^1}$ we have 
\begin{align}\label{Cond-Tiempo}
		T\big(\norm{w_{0}}_{L^1}\big) \geq T_1. 
\end{align}
Then,  for $0 < \varepsilon < T_1$ small enough, we consider the time $T^* - \varepsilon > 0$ and we set the initial datum  $w_0 = u_\alpha(T^* - \varepsilon, \cdot)$. We  shall denote by $w_\alpha$ its  corresponding arising solution, which exists at least until the time $T\big(\norm{w_{0}}_{L^1}\big)$. We thus observe that the function 
 \begin{align*}
 		\tilde{u}_\alpha(t, \cdot)= \left\{\begin{array}{l l}
 			\displaystyle u_\alpha(t, \cdot),
 			\qquad  t \in \big[0, T^* - \varepsilon \big],\\ \\
 			\displaystyle w_\alpha(t, \cdot),
 			\qquad  t \in \big[T^* - \varepsilon, T^* - \varepsilon+T\big(\norm{v_{0}}_{L^1}\big)\big],
 		\end{array}\right.
 	\end{align*}
is a solution of the equation \eqref{u-alpha} associated to the initial datum $u_{0}$. Moreover, we observe that this function is defined in the interval of time  $\big[0, T^* - \varepsilon+T\big(\norm{v_0}_{L^1}\big)\big]$.  But, by \eqref{eq:15} we have $\ds{	\norm{w_{0}}_{L^1} = \norm{u_\alpha(T^* - \varepsilon, \cdot)}_{L^1} \leq \norm{u_0}_{L^1}}$ and consequently  by (\ref{Cond-Tiempo}) we get  $T\big(\norm{w_{0}}_{L^1}\big) \geq T_1$. 

\medskip

Finally, we can write $	T^* - \varepsilon + T_1 \leq T^* - \varepsilon + T\big(\norm{v_{0}}_{L^1}\big)$;  and since $0 < \varepsilon < T_1$  we obtain $T^*<T^* - \varepsilon + T_1$, which is  a contradiction with the definition of the time $T^*$.  Proposition  \ref{Th-Tech-4} is proven. \finpv

\medskip 

Once we have proven the Propositios \ref{Th-tech-1}, \ref{prop:6} and \ref{Th-Tech-4}, we  can finish with the proof of Theorem \ref{Th-WP} . \finpv

\section{From anomalous to classical diffusion}\label{Sec:Anom-to-local}
\subsection{Proof of  Theorem \ref{Th-Unif-Conv}}
For $1<\alpha <2$, let $u_\alpha$ be the mild solution of the equation (\ref{Eq3})  given by the expression (\ref{u-alpha}). Moreover, for $\alpha=2$ let $u_2$ be the mild solution of the equation (\ref{Eq3}), which is given by the expression (\ref{u-2}). Then, for a time $0<T<+\infty$ fixed  we write 
\begin{equation}\label{estim-base}
\begin{split}
\sup_{0\leq t \leq T}&\Vert u_{\alpha}(t,\cdot)-u_2(t,\cdot)\Vert_{L^{\infty}}\leq  \sup_{0\leq t \leq T} \left\Vert p_\alpha(t,\cdot)\ast u_{0,\alpha} -h(t,\cdot)\ast u_{0,2} \right\Vert_{L^{\infty}} \\
&+  \sup_{0\leq t \leq T}\left\Vert \int_{0}^{t}p_\alpha(t-s,\cdot)\ast \eta \cdot \nabla(u^{b}_{\alpha})(s,\cdot) ds- \int_{0}^{t}h(t-s,\cdot)\ast \eta \cdot \nabla(u^{b}_{2})(s,\cdot) ds \right\Vert_{L^{\infty}}\\
&=I_{\alpha}+J_{\alpha},
\end{split}
\end{equation}
where we must estimate each term on the right-hand side. For the  term $I_\alpha$   we write 
\begin{equation}\label{estim-lin-01}
I_\alpha \leq  \sup_{0\leq t \leq T} \left\Vert \left( p_\alpha(t,\cdot) - h(t,\cdot) \right)\ast u_{0,\alpha} \right\Vert_{L^\infty} +  \sup_{0\leq t \leq T} \left\Vert h(t,\cdot)\ast \left( u_{0,\alpha} - u_{0,2}\right) \right\Vert_{L^\infty} = I_{\alpha,1}+I_{\alpha,2}.     
\end{equation}
In order to estimate the term $I_{\alpha,1}$, we apply  the Bessel potential operators $(1-\Delta)^{-s/2}$ and $(1-\Delta)^{s/2}$ to obtain:  
\begin{equation*}
I_{\alpha,1} = \sup_{0\leq t \leq T}\left\Vert (1-\Delta)^{-s/2}\Big(p_\alpha(t,\cdot)-h(t,\cdot)\Big) \ast (1-\Delta)^{s/2} u_{0,\alpha} \right\Vert_{L^{\infty}}=(a).   
\end{equation*}
Then, by applying the  Young inequalities (with $1+1/\infty=1/2+1/2$) we have 
\begin{equation}\label{estim01}
\begin{split}
(a) &\leq  c\,  \sup_{0\leq t \leq T} \left(\left\Vert (1-\Delta)^{-s/2}\Big(p_\alpha(t,\cdot)-h(t,\cdot)\Big)  \right\Vert_{L^{2}}\, \left\Vert (1-\Delta)^{s/2} u_{0,\alpha} \right\Vert_{L^{2}}\right) \\
& \leq  c \left(  \sup_{0\leq t \leq T} \left\Vert p_\alpha(t,\cdot)-h(t,\cdot) \right\Vert_{H^{-s}} \right)\, \left(  \sup_{1 < \alpha < 2 } \Vert u_{0,\alpha} \Vert_{H^{s}} \right),
\end{split}
\end{equation}
where we shall control each term  above separately. For the first term on the right-hand side we have the following technical result:
\begin{Lemme}\label{Lema-Tech-2} For $s>n/2$   there exists a constant $C=C(s)>0$  such that  for all $1<\alpha < 2$ we  have:  
\[ \sup_{0\leq t \leq T} \left\Vert p_\alpha(t,\cdot)-h(t,\cdot) \right\Vert_{H^{-s}}\leq  C \, T\, \vert 2- \alpha \vert.\]
\end{Lemme}
\pv First,  we verify that the quantity $\ds{ \left\Vert p_\alpha(t,\cdot)-h(t,\cdot) \right\Vert^{2}_{H^{-s}}}$ is continuous in the temporal variable $t$. Indeed, for  $0\leq t_0, t\leq T$ we have 
\begin{equation*} 
\begin{split}
 & \left\Vert p_\alpha(t,\cdot)-h(t,\cdot) \right\Vert^{2}_{H^{-s}}   -   \left\Vert p_\alpha(t_0,\cdot)-h(t_0,\cdot) \right\Vert^{2}_{H^{-s}}  \\
 =& \int_{\Rn} \left\vert  e^{-\vert \xi \vert^{\alpha} t} -e^{-\vert \xi \vert^2 t}\right\vert^2 \frac{d \xi}{(1+\vert \xi \vert^2)^{s}} - \int_{\Rn} \left\vert  e^{-\vert \xi \vert^{\alpha} t_0} -e^{-\vert \xi \vert^2 t_0}\right\vert^2 \frac{d \xi}{(1+\vert \xi \vert^2)^{s}}\\ =&  \int_{\Rn} \left( \left\vert  e^{-\vert \xi \vert^{\alpha} t} -e^{-\vert \xi \vert^2 t}\right\vert^2- \left\vert  e^{-\vert \xi \vert^{\alpha} t_0} -e^{-\vert \xi \vert^2 t_0}\right\vert^2 \right)\frac{d \xi}{(1+\vert \xi \vert^2)^{s}}.
\end{split} 
\end{equation*}
As  $s>n/2$  we have $\ds{\int_{\Rn}\frac{d \xi}{(1+\vert \xi \vert^2)^{s}} <+\infty}$; and then,  we can apply the dominated  convergence theorem  to obtain that $\ds{\lim_{t \to t_0} \left( \left\Vert p_\alpha(t,\cdot)-h(t,\cdot) \right\Vert^{2}_{H^{-s}}   -   \left\Vert p_\alpha(t_0,\cdot)-h(t_0,\cdot) \right\Vert^{2}_{H^{-s}} \right)=0}$. 

\medskip

Thereafter, by the continuity of the quantity $\ds{ \left\Vert p_\alpha(t,\cdot)-h(t,\cdot) \right\Vert^{2}_{H^{-s}}}$ with respect to the variable $t$,   there exists a time $0 < t_1 \leq T$ such that $\ds{\sup_{0\leq t \leq T} \left\Vert p_\alpha(t,\cdot)-h(t,\cdot) \right\Vert_{H^{-s}}=\left\Vert p_\alpha(t_1,\cdot)-h(t_1,\cdot) \right\Vert_{H^{-s}}}$. 
\medskip
Now, we will prove the estimate $\ds{\left\Vert p_\alpha(t_1,\cdot)-h(t_1,\cdot) \right\Vert_{H^{-s}}\leq  C\, T\,  \vert 2- \alpha \vert}$. For this we  write: 
\begin{equation}\label{Id}
 \left\Vert p_\alpha(t_1,\cdot)-h(t_1,\cdot) \right\Vert^{2}_{H^{-s}} = \int_{\Rn} \vert  e^{-\vert \xi \vert^{\alpha} t_1} -e^{-\vert \xi \vert^2 t_1}\vert^2 \frac{d \xi}{(1+\vert \xi \vert^2)^{s}}.
\end{equation} 
Here, for   $\xi \in \Rn \setminus \{0\}$ fixed, and for $1<\alpha < 2 +\delta$ (with $\delta>0$) we define the function 
\begin{equation}\label{function}
 f_\xi (\alpha)= e^{- t_1 \vert \xi \vert^\alpha},   
\end{equation}
where, by computing its derivative with respect to the variable $\alpha$ we get
\[ f^{'}_{\xi}(\alpha)=- t_1\,  e^{- t_1 \vert \xi \vert^\alpha} \, \vert \xi \vert^\alpha \ln (\vert \xi \vert).\] 
Thus, by the mean value theorem (in the variable $\alpha$)   we can write 
\[ \vert f_\xi(\alpha)-f_\xi(2) \vert \leq \Vert f^{'}_{\xi} \Vert_{L^{\infty}([1,2+\delta])}\, \vert 2-\alpha \vert.\] 
Moreover,  we can  also prove the uniform estimate with respect to the variable $\xi$:
\begin{equation}\label{Estim-Tech}
\left\Vert  \Vert f^{'}_{\xi} \Vert_{L^{\infty}([1,2+\delta])} \right\Vert_{L^{\infty}(\Rn)}  \leq c\, T.    
\end{equation}
The proof of this estimate is not difficult and it is given in detail at the Appendix \ref{AppendixB}. We thus have, 
\[ \vert f_\xi(\alpha) - f_\xi(2) \vert \leq c\, T \vert 2-\alpha \vert. \]
Then, we get  back to the identity (\ref{Id}) and we can write 
\begin{equation*}
\Vert p_\alpha(t_1, \cdot) - h(t_1,\cdot) \Vert^{2}_{H^{-s}} = \int_{\Rn} \vert f_\xi (\alpha) - f_\xi(2) \vert^2 \frac{d \xi}{(1+\vert \xi \vert^2)^s} \leq c\, T^2\, \vert 2-\alpha \vert^2\,  \int_{\Rn}  \frac{d \xi}{(1+\vert \xi \vert^2)^s} \leq \, C(s) \, T^2\, \vert 2-\alpha \vert^2.
\end{equation*}
\finpv
 
\medskip

On the other hand, to control the second term on the right-hand side in the expression (\ref{estim01}), we ust recall that by hypothesis (\ref{Conv-Data})  the family $\ds{(u_{0,\alpha})_{1<\alpha <2}}$ is bounded in $H^s(\Rn)$.  

\medskip 

Thus, for the term $I_{\alpha,1}$ given in  (\ref{estim-lin-01})  we can write:  
\begin{equation}\label{Lim1-01}
 I_{\alpha,1} \leq C\, T\,  \vert 2-\alpha\vert.    
\end{equation}

We study now the term $I_{\alpha,2}$, which is also given in (\ref{estim-lin-01}).  By the Young inequalities (with $1+ 1/\infty= 1 + 1/ \infty)$, the well-known properties of the heat kernel, and moreover, by the assumption give in (\ref{Conv-Rate-Data}),  we have: 
\begin{equation}\label{Lim-02}
I_{\alpha,2} \leq {\bf c} (2-\alpha)^\gamma.     
\end{equation}
Consequently, with the  estimates  (\ref{Lim1-01}) and (\ref{Lim-02}) above we obtain:
\begin{equation}\label{Lin}
I_\alpha \leq {\bf C}\,(1+T)\, \max\Big( (2-\alpha)^\gamma, 2-\alpha \Big).
\end{equation}

We study now the term $J_{\alpha}$ given in the expression (\ref{estim-base}). For this we write 
\begin{equation}\label{Estim-J}
\begin{split}
J_{\alpha} \leq &  \sup_{0\leq t \leq T}  \left\Vert \int_{0}^{t}p_\alpha(t-s,\cdot)\ast \eta \cdot \nabla(u^{b}_{\alpha})(s,\cdot) ds- \int_{0}^{t}h_\alpha(t-s,\cdot)\ast \eta \cdot  \nabla(u^{b}_{\alpha})(s,\cdot) ds\right\Vert_{L^{\infty}} \\
&+ \sup_{0\leq t \leq T}  \left\Vert \int_{0}^{t}h(t-s,\cdot)\ast \eta \cdot \nabla(u^{b}_{\alpha})(s,\cdot) ds- \int_{0}^{t}h(t-s,\cdot)\ast \eta \cdot \nabla(u^{b}_{2})(s,\cdot) ds\right\Vert_{L^{\infty}} \\
\leq & \sup_{0\leq t \leq T}  \left\Vert \int_{0}^{t}\Big(p_\alpha(t-s,\cdot)-h(t-s,\cdot)\Big)\ast \eta \cdot \nabla(u^{b}_{\alpha})(s,\cdot) ds\right\Vert_{L^{\infty}} \\
&+ \sup_{0\leq t \leq T} \left\Vert \int_{0}^{t}  h(t-s,\cdot)\ast \eta \cdot \nabla \Big(u^{b}_{\alpha}-u^{b}_{2} \Big) (s,\cdot) ds \right\Vert_{L^{\infty}}=J_{\alpha,1}+J_{\alpha,2},
\end{split}    
\end{equation} 
where  we will  study the terms $J_{\alpha,1}$ and $J_{\alpha,2}$ separately. For the term $J_{\alpha,1}$, we apply  first the operators $(1-\Delta)^{-s/2}$ and $(1-\Delta)^{s/2}$, and moreover,  by  the Young inequalities (with $1+1/\infty= 1/2 + 1/2$)  we have 
\begin{equation}\label{estim-J1}  
\begin{split}
J_{\alpha,1} & \leq  \sup_{0\leq t \leq T} \left( \int_{0}^{t} \left\Vert \Big(p_\alpha(t-s,\cdot)-h(t-s,\cdot)\Big)\ast \eta \cdot \nabla (u^{b}_{\alpha})(s,\cdot) \right\Vert_{L^{\infty}} ds   \right) \\
&\leq \vert \eta \vert \sup_{0\leq t \leq T} \left( \int_{0}^{t} \left\Vert \nabla p_\alpha(t-s,\cdot)-\nabla h(t-s,\cdot)\right\Vert_{H^{-s}}\, \left\Vert  u^{b}_{\alpha}(s,\cdot) \right\Vert_{H^s} ds  \right) \\
&\leq \vert \eta \vert\,  T \left( \sup_{0\leq t \leq T} \left\Vert \nabla  p_\alpha(t,\cdot)- \nabla h(t,\cdot)\right\Vert_{H^{-s}}\right) \left( \sup_{0\leq t \leq T} \Vert u^{b}_{\alpha}(s,\cdot) \Vert_{H^s}\right).
\end{split}    
\end{equation}

In order to control the first term on the right-hand side, we follow the same lines in the proof of Lemma \ref{Lema-Tech-2} with the function $\ds{f_\xi(\alpha)= i \xi_j \,  e^{- t_1 \vert \xi \vert^\alpha}}$, with $j=1,2,\cdots, n$. Then we have 
\begin{equation}\label{Estim-grad-kernels}
\sup_{0\leq t \leq T} \left\Vert \nabla p_\alpha(t,\cdot)- \nabla h(t,\cdot)\right\Vert_{H^{-s}} \leq C\, T\,  \vert 2-\alpha \vert.
\end{equation}

Thereafter, to control the second term on the right-hand side we shall need the following: 
\begin{Lemme}\label{Lema-Tech-1} There exists $0<\varepsilon\ll 1$, and there exists a constant $\mathcal{C}=\mathcal{C}(\varepsilon,T,b, \Vert u_{0,2}\Vert_{L^1}, \Vert u_{0,2}\Vert_{H^s})>0$, such that for all $1+\varepsilon < \alpha < 2$ we have:  
\begin{equation}\label{Estim-unif-alpha}
\sup_{0\leq t \leq T} \Vert u^{b}_{\alpha}(t,\cdot) \Vert_{H^s} \leq \mathcal{C}.
\end{equation}	
\end{Lemme}
\pv For any the initial data $u_{0,\alpha} \in L^1 \cap H^s(\Rn)$, with $1<\alpha<2$,  we recall that by Proposition \ref{Th-tech-1} there exists a  time $T_\alpha$ (depending on $\alpha$)  defined by (\ref{Time-T}) as 
\[ T_\alpha =  \dfrac{1}{2} \left[\dfrac{1-\frac{1}{\alpha}}{2^b c \abs{\eta} \Big(\norm{u_{0,\alpha}}_{L^1} + \norm{u_{0,\alpha}}_{H^s} \Big)^{b-1}}\right]^{\frac{\alpha}{\alpha-1}}, \]
and there exists  a (unique) solution $u_\alpha \in E_{T_\alpha}$ of the equation (\ref{u-alpha}). Our staring point is to obtain a lower bound for the time $T_\alpha$ which does not depend on $\alpha$. 

\medskip

By our hypothesis given in  (\ref{Conv-Data}) we can set $0<\varepsilon \ll 1$ such that for all $1+\varepsilon < \alpha <2$ we have: 
\[ \left\vert \Big( \Vert u_{0,\alpha}\Vert_{L^1}+\Vert u_{0,\alpha}\Vert_{H^s}\Big) - \Big( \Vert u_{0,2}\Vert_{L^1}+\Vert u_{0,2}\Vert_{H^s} \Big)   \right\vert \leq \frac{1}{2} \Big( \Vert u_{0,2}\Vert_{L^1}+\Vert u_{0,2}\Vert_{H^s} \Big),\] 
hence we get the control 
\begin{equation}\label{Control-alpha-2}
\Big( \Vert u_{0,\alpha}\Vert_{L^1}+\Vert u_{0,\alpha}\Vert_{H^s} \Big) \leq  \frac{3}{2} \Big( \Vert u_{0,2}\Vert_{L^1}+\Vert u_{0,2}\Vert_{H^s}\Big), \quad 1+\varepsilon < \alpha < 2, 
\end{equation}
and then, we can write: 
\begin{equation*}
 \frac{1}{2}\left[ \frac{1-\nicefrac{1}{\alpha}}{2^b \, c \vert \eta \vert (\Vert u_{0,2} \Vert_{L^1}+\Vert u_{0,2} \Vert_{H^s})^{b-1}} \right]^{\nicefrac{\alpha}{\alpha-1}}  \leq T_\alpha.
\end{equation*}
Moreover, since  $1+\varepsilon < \alpha <2$  then the expression on the left-hand  can be estimated from below  by the following quantity: 
\begin{equation}\label{T0}
T_0=\max\left( \frac{1}{2}\left[ \frac{1-\nicefrac{1}{1+\varepsilon}}{2^b \, c \vert \eta \vert (\Vert u_{0,2} \Vert_{L^1}+\Vert u_{0,2} \Vert_{H^s})^{b-1}} \right]^{2/\varepsilon},  \frac{1}{2}\left[ \frac{1-\nicefrac{1}{1+\varepsilon}}{2^b \, c \vert \eta \vert (\Vert u_{0,2} \Vert_{L^1}+\Vert u_{0,2} \Vert_{H^s})^{b-1}} \right]^{1+\varepsilon} \right).
\end{equation}
At the Appendix \ref{AppendixA}  we verify in detail this estimate. Then,  for all $1+\varepsilon<\alpha < 2$ we have $T_0  \leq T_\alpha$. 

\medskip

Once we have the lower estimate $T_0\leq T_\alpha$, we remark that for all $1+\varepsilon < \alpha < 2$ the solution $u_\alpha$ of the equation (\ref{u-alpha}), which is  constructed in the Proposition \ref{Th-tech-1} by the Picard's fixed point argument, verifies $u_\alpha \in E_{T_0}$  and moreover we have the estimate  $\ds{\Vert u_\alpha \Vert_{E_{T_0}} \leq c_0(\Vert u_{0,\alpha} \Vert_{L^1}+\Vert u_{0,\alpha} \Vert_{H^s})}$.

\medskip

We also remark that by Proposition \ref{Th-Tech-4} the solution $u_\alpha$ is extended to a global in time solution  by a well-known iterative argument: for every interval $[kT_0, (k+1)T_0]$  (with $k \in \mathbb{N}^{*}$)  we set the initial initial datum $u_\alpha(k T_0, \cdot)$ and we apply again the Picard's fixed point schema to obtain a (unique) solution  $u_\alpha \in E_{[k T_0, (k+1)T_0]}$ (recall that the space  $ E_{[k T_0, (k+1)T_0]}$ is defined in  (\ref{Espace}) and (\ref{Norm})).  Moreover,  there  exists a constant $c_k >0$ such that we have  $\ds{\Vert u_\alpha \Vert_{E_{[kT_0, (k+1)T_0]}} \leq c_k (\Vert u_\alpha(k T_0, \cdot) \Vert_{L^1}+\Vert u_\alpha(k T_0, \cdot) \Vert_{H^s})}$. 

\medskip 

We study now the expression $\ds{c_k (\Vert u_\alpha(k T_0, \cdot) \Vert_{L^1}+\Vert u_\alpha(k T_0, \cdot) \Vert_{H^s})}$.  For the quantity $\ds{\Vert u_\alpha(k T_0, \cdot) \Vert_{L^1} }$, by (\ref{eq:15}) we have $\ds{\Vert u_\alpha (kT_0, \cdot)\Vert_{L^1} \leq \Vert u_{0,\alpha} \Vert_{L^1}}$. Then, we can write 
\[ c_k (\Vert u_\alpha(k T_0, \cdot) \Vert_{L^1}+\Vert u_\alpha(k T_0, \cdot) \Vert_{H^s}) \leq c_k (\Vert u_{0,\alpha} \Vert_{L^1}+\Vert u_\alpha(k T_0, \cdot) \Vert_{H^s}).\]
On the other hand, for the quantity $\ds{\Vert u_\alpha(k T_0, \cdot  ) \Vert_{H^s}}$, we remark that we have 
\[  \Vert u_\alpha(k T_0, \cdot  ) \Vert_{H^s} \leq \sup_{(k-1)T_0 \leq t \leq k T_0} \Vert  u_\alpha(t,\cdot) \Vert_{H^s} \leq \Vert u_\alpha \Vert_{E_{[(k-1)T_0, k T_0]}} \leq c_{k-1} (\Vert u_{0,\alpha} \Vert_{L^1} + \Vert u_{\alpha}((k-1)T_0, \cdot)) \Vert_{H^s}.\] 
Thus, we can iterate these estimates and for a  constant $C_k>0$ big enough (in particular we must have  $\ds{C_k > \prod_{j=0}^{k} c_j}$) we obtain  $\Vert u_\alpha(k T_0, \cdot  ) \Vert_{H^s} \leq C_k (\Vert u_{0,\alpha} \Vert_{L^1}+\Vert u_{0,\alpha} \Vert_{H^s})$.  Consequently,  for all $k\in \mathbb{N}^{*}$ we have 
\begin{equation}\label{Control-initial-data}
\Vert u_\alpha \Vert_{E_{[k T_0, (k+1)T_0]}} \leq C_k (\Vert u_{0,\alpha} \Vert_{L^1}+\Vert u_{0,\alpha} \Vert_{H^s}).
\end{equation}

Now, we are able to prove the estimate (\ref{Estim-unif-alpha}) stated in this lemma. For the  time $T$  there exists $k_T \in \mathbb{N}$ (which depends on $T$) such that we have $k_T T_0 \leq T \leq (k_T+1) T_0$. Then,  as $s>n/2$ by the product laws in the Sobolev spaces we can write 
\begin{equation*}
\begin{split}
 \sup_{0 \leq t \leq T} \Vert u^{b}_\alpha(t,\cdot)\Vert_{H^s}  \leq &  \sup_{0 \leq t \leq T} \Vert u_\alpha(t,\cdot)\Vert^{b}_{H^s} 
  \leq \, \left( \sup_{0 \leq t \leq T} \Vert u_\alpha(t,\cdot)\Vert_{H^s}  \right)^{b}  \\
  \leq &\, \left( \sum_{j=0}^{k_T} \, \sup_{j T_0 \leq (j+1) T_0} \Vert u_\alpha (t,\cdot)\Vert_{H^s} \right)^b  \leq   \left( \sum_{j=0}^{k_T} \,  \Vert u_\alpha \Vert_{E_{[jT_0, (j+1)T_0}} \right)^b.    
\end{split}      
\end{equation*}
Then, by the control given in (\ref{Control-initial-data})  and the control given in (\ref{Control-alpha-2}) we have 
\begin{equation*}
\begin{split}
\left( \sum_{j=0}^{k_T} \,  \Vert u_\alpha \Vert_{E_{[jT_0, (j+1)T_0}} \right)^b \leq  &\,  \left( \sum_{j=0}^{k_T} \,  C_j (\Vert u_{0,\alpha} \Vert_{L^1}+ \Vert u_{0,\alpha} \Vert_{H^s}) \right)^b 
\leq  \left( \sum_{j=0}^{k_T}\, C_j  \right)^b \, (\Vert u_{0,\alpha} \Vert_{L^1}+ \Vert u_{0,\alpha} \Vert_{H^s})^b \\
\leq &\,  \left( \sum_{j=0}^{k_T}\, C_j  \right)^b\, c\, (\Vert u_{0,2}\Vert_{L^1}+\Vert u_{0,2} \Vert_{H^s})^b
=  \mathcal{C}(\varepsilon, b,\Vert u_{0,2}\Vert_{L^1},  \Vert u_{0,2} \Vert_{H^s}, T).    
\end{split}      
\end{equation*}
To finish the proof of this lemma, we just remark that the constant $\mathcal{C}$ defined above also depends on the parameter $\varepsilon$, since the $k_T$ depends on $T$ and $T_0$; and the time $T$  given in (\ref{T0}) depends on $\varepsilon$. \finpv
 
\medskip 

With the estimates (\ref{Estim-grad-kernels}) and (\ref{Estim-unif-alpha}) at our disposal,  we get back to the estimate (\ref{estim-J1}) to  write 
\begin{equation}\label{Lim2}
 J_{\alpha,1} \leq C\, \vert \eta\vert \, T^2 \, \vert 2-\alpha \vert \leq  {\bf C}\, \vert \eta\vert \, T ^2\, \max\Big( (2-\alpha)^\gamma, 2-\alpha \Big).    
\end{equation}

On the other hand, We study now  the term $J_{\alpha,2}$ given in (\ref{Estim-J}). For this, by  the Young inequalities  we write:  
\[ J_{\alpha,2} \leq c\, \vert \eta \vert \, \sup_{0\leq t \leq T} \int_{0}^{t} \Vert \nabla h(t-s,\cdot)\Vert_{L^1} \Vert u^{b}_{\alpha}(s,\cdot)- u^{b}_{2}(s,\cdot) \Vert_{L^\infty} ds=(a). \]
Here, by the well-known properties of the heat kernel $h(t,\cdot)$ we have $\ds{\Vert \nabla h(t-s,\cdot)\Vert_{L^1} \leq c (t-s)^{-1/2}}$. Thereafter,  in order to estimate the term $\ds{\Vert u^{b}_{\alpha}(s,\cdot)- u^{b}_{2}(s,\cdot) \Vert_{L^\infty}}$, since  $s>n/2$ and  by Lemma \ref{Lema-Tech-2}  we can write 
\begin{equation*}
\begin{split}
\Vert u^{b}_{\alpha}(s,\cdot)- u^{b}_{2}(s,\cdot) \Vert_{L^\infty} =& \left\Vert (u_{\alpha}(s,\cdot)-u_{2}(s,\cdot)) \, \sum_{j=0}^{b-1} u^{b-1-j}_{\alpha}(s,\cdot) u^{j}_{2}(s,\cdot) \right\Vert_{L^\infty} \\
\leq &\Vert u_{\alpha}(s,\cdot)-u_{2}(s,\cdot) \Vert_{L^\infty} \,  \sum_{j=0}^{b-1}\left\Vert u_{\alpha}(s,\cdot)\right\Vert^{b-1-j}_{L^\infty} \, \left\Vert u_{2}(s,\cdot) \right\Vert^{j}_{L^\infty} \\
\leq &\Vert u_{\alpha}(s,\cdot)-u_{2}(s,\cdot) \Vert_{L^\infty} \,   \sum_{j=0}^{b-1}\left\Vert u_{\alpha}(s,\cdot)\right\Vert^{b-1-j}_{H^s} \, \left\Vert u_{2}(s,\cdot) \right\Vert^{j}_{H^s}\\
\leq &  C\,  \Vert u_{\alpha}(s,\cdot)-u_{2}(s,\cdot) \Vert_{L^\infty},
\end{split}    
\end{equation*}
where the constant $C>0$ does not depend on $\alpha$.  With these estimates we obtain: 
\begin{equation*}
\begin{split}
(a) \leq & C\, \vert \eta \vert \, \sup_{0\leq t \leq T} \, \int_{0}^{t}(t-s)^{-1/2} \, \Vert  u_{\alpha}(s,\cdot)-u_{2}(s,\cdot) \Vert_{L^\infty} ds \leq C\, \vert \eta \vert T^{1/2}\, \left( \sup_{0\leq s \leq T} \Vert  u_{\alpha}(s,\cdot)-u_{2}(s,\cdot) \Vert_{L^\infty}  \right), 
\end{split}    
\end{equation*}
and consequently we have 
\begin{equation}\label{Lim3}
J_{\alpha,2} \leq   C\, \vert \eta \vert T^{1/2}\, \left( \sup_{0\leq s \leq T} \Vert  u_{\alpha}(s,\cdot)-u_{2}(s,\cdot) \Vert_{L^\infty}  \right).
\end{equation}

Once we estimated  the terms $I_\alpha$, $J_{\alpha,1}$ and $J_{\alpha,2}$ in (\ref{Lin}), (\ref{Lim2}) and (\ref{Lim3}) respectively, we  get back to (\ref{estim-base}) we can write 
\[  \sup_{0\leq t \leq T} \Vert u_\alpha(t,\cdot)- u_2(t,\cdot)\Vert_{L^\infty} \leq I_\alpha + J_{\alpha,1}+J_{\alpha, 2}  \leq I_\alpha + J_{\alpha,1}+  C\, \vert \eta \vert \,  T^{1/2}\, \left( \sup_{0\leq s \leq T} \Vert  u_{\alpha}(s,\cdot)-u_{2}(s,\cdot) \Vert_{L^\infty}  \right). \]

In this estimate, first we set a  time $0<T_1<T$ small enough such that it verifies  
\begin{equation*}
C\, \vert \eta \vert \,  T^{1/2}_{1} \leq \frac{1}{2},   
\end{equation*}
we get: 
\[  \sup_{0\leq t \leq T_1} \Vert u_\alpha(t,\cdot)- u_2(t,\cdot)\Vert_{L^\infty}  \leq I_\alpha + J_{\alpha,1}+  \frac{1}{2}\, \left( \sup_{0\leq s \leq T_1} \Vert  u_{\alpha}(s,\cdot)-u_{2}(s,\cdot) \Vert_{L^\infty}  \right), \]
and then we can write 
\[ \frac{1}{2} \sup_{0\leq t \leq T_1} \Vert u_\alpha(t,\cdot)- u_2(t,\cdot)\Vert_{L^\infty} \leq I_\alpha + J_{\alpha,1}.\]
Hence, by (\ref{Lin}) and (\ref{Lim2}) we obtain
\begin{equation*}
\begin{split}
\sup_{0\leq t \leq T_1} \Vert u_\alpha(t,\cdot)- u_2(t,\cdot)\Vert_{L^\infty} \leq  & \, {\bf C}(1+T_1+T^{2}_{1}) \, \max\left( (2-\alpha)^\gamma, (2-\alpha)\right) \\
\leq & \,  {\bf C}(1+T+T^{2}) \, \max\left( (2-\alpha)^\gamma, (2-\alpha)\right). 
\end{split}
\end{equation*}

Finally, we  iterate this argument on the intervals $[kT_1, (k+1)T_1]$, with $k \in \mathbb{N}$, and then,   for  the time $0<T<+\infty$ we have 
\[ \sup_{0\leq t \leq T} \Vert u_\alpha(t,\cdot)- u_2(t,\cdot)\Vert_{L^\infty} \leq  {\bf C}\, (1+T+T^2)  \, \max \Big( (2-\alpha)^\gamma, 2-\alpha \Big ).\]
Theorem  \ref{Th-Unif-Conv} is now proven. \finpv   

\subsection{Proof of the Corollary \ref{Col-1}} 
For $0<T<+\infty$ fixed, and moreover, for  $1\leq q < +\infty$ and $1< p < +\infty$, by the interpolation inequalities (with $\theta = 1/q$) we write 
\[ \left( \int_{0}^{T} \Vert u_\alpha(t,\cdot)-u_2(t,\cdot)\Vert^{q}_{L^p} dt \right)^{1/q} \leq \left( \int_{0}^{T} \Vert u_\alpha(t,\cdot)-u_2(t,\cdot)\Vert^{q\, \theta}_{L^1}\, \Vert u_\alpha(t,\cdot)-u_2(t,\cdot)\Vert^{q (1-\theta)}_{L^\infty}\, dt \right)^{1/q}=(a).\]
By the hypothesis (\ref{Conv-Data}), we can set  $M>0$ (small enough) such that for $1+\varepsilon < \alpha <2$  we have  $\ds{\Vert u_\alpha(t,\cdot)-u_2(t,\cdot)\Vert_{L^1} \leq M}$. Then  we obtain 
\[ (a) \leq M^{\theta}\, \left( \int_{0}^{T} \Vert u_\alpha(t,\cdot)-u_2(t,\cdot)\Vert^{q(1-\theta)}_{L^\infty} dt \right)^{1/q} \leq M^{\theta} \Vert u_\alpha(t,\cdot)-u_2(t,\cdot)\Vert^{(1-\theta)}_{L^\infty}\, T^{1/q}.\]
Thus, the wished estimate follows from the inequality (\ref{Conv-Rate-Sol}) proven in the Theorem \ref{Th-Unif-Conv}.  \finpv

\begin{appendices}
\section{Appendix}\label{AppendixB}
We prove here the estimate (\ref{Estim-Tech}). We recall the expression 
\[ f^{'}_{\xi}(\alpha)=-t_1 e^{-t_1 \vert \xi \vert^{\alpha}} \vert \xi \vert^\alpha \ln (\vert \xi \vert), \quad 1<\alpha <2+\delta, \quad 0<t_1\leq T.\]
Then, we write  
\[ \left\Vert  \Vert f^{'}_{\xi} \Vert_{L^{\infty}([1,2+\delta])} \right\Vert_{L^{\infty}(\Rn)}\leq \left\Vert  \Vert f^{'}_{\xi} \Vert_{L^{\infty}([1,2+\delta])} \right\Vert_{L^{\infty}(\vert \xi \vert \leq 1)}+ \left\Vert  \Vert f^{'}_{\xi} \Vert_{L^{\infty}([1,2+\delta])} \right\Vert_{L^{\infty}(\vert \xi \vert >1)}=A+B,\]
where, we shall estimate the terms $A$ and $B$ separately. For the term $A$,  as we have $\vert \xi \vert \leq 1$, $1<\alpha<2+\delta$, and moreover, as we have $\ds{\lim_{\vert \xi \vert\to 0^{+}} \vert \xi \vert \ln(\vert \xi \vert)=0}$, then we can write:
\[ A \leq  T \, \left( \sup_{\xi \in \Rn} e^{-t_1 \vert \xi \vert^{2+\delta}} \vert \xi \vert \ln (\vert \xi \vert)\right)\leq C\, T.   \]
For the term $B$, since $\vert \xi \vert >1$ then we can write 
\[B \leq T\, \left( \sup_{\xi \in \Rn} e^{-t_1 \vert \xi \vert} \vert \xi \vert^{2+\delta} \ln (\vert \xi \vert)\right)\leq C\, T. \]

\section{Appendix}\label{AppendixA}
Here we give a proof of the estimate 
\begin{equation*}
\begin{split}
 &T_0=\max\left( \frac{1}{2}\left[ \frac{1-\nicefrac{1}{1+\varepsilon}}{2^b \, c \vert \eta \vert (\Vert u_{0,2} \Vert_{L^1}+\Vert u_{0,2} \Vert_{H^s})^{b-1}} \right]^{2/\varepsilon},  \frac{1}{2}\left[ \frac{1-\nicefrac{1}{1+\varepsilon}}{2^b \, c \vert \eta \vert (\Vert u_{0,2} \Vert_{L^1}+\Vert u_{0,2} \Vert_{H^s})^{b-1}} \right]^{1+\varepsilon} \right)\\
 \leq &  \frac{1}{2}\left[ \frac{1-\nicefrac{1}{\alpha}}{2^b \, c \vert \eta \vert (\Vert u_{0,2} \Vert_{L^1}+\Vert u_{0,2} \Vert_{H^s})^{b-1}} \right]^{\nicefrac{\alpha}{\alpha-1}}.  
\end{split}    
\end{equation*}
First, as we have $1+\varepsilon < \alpha < 2$, then we get $1-\frac{1}{1+\varepsilon}< 1-\frac{1}{\alpha}$, and we can write 
\[ \frac{1}{2} \left[ \frac{1-\nicefrac{1}{1+\varepsilon}}{2^b c \vert \eta \vert (\Vert u_{0,2} \Vert_{L^1}+\Vert u_{0,2} \Vert_{H^s})^{b-1}}\right]^{\frac{\alpha}{\alpha-1}} \leq \frac{1}{2}\left[ \frac{1-\nicefrac{1}{\alpha}}{2^b \, c \vert \eta \vert (\Vert u_{0,2} \Vert_{L^1}+\Vert u_{0,2} \Vert_{H^s})^{b-1}} \right]^{\nicefrac{\alpha}{\alpha-1}}.\]
Thereafter, by the sake of simplicity, we denote 
\[  \frac{1-\nicefrac{1}{1+\varepsilon}}{2^b c \vert \eta \vert (\Vert u_{0,2} \Vert_{L^1}+\Vert u_{0,2} \Vert_{H^s})^{b-1}}= (a),\]
and we have 
\[ \frac{1}{2} [(a)]^{\frac{\alpha}{\alpha-1}} \leq \frac{1}{2}\left[ \frac{1-\nicefrac{1}{\alpha}}{2^b \, c \vert \eta \vert (\Vert u_{0,2} \Vert_{L^1}+\Vert u_{0,2} \Vert_{H^s})^{b-1}} \right]^{\nicefrac{\alpha}{\alpha-1}}.\]
We study now the expression $\frac{\alpha}{\alpha-1}$.  Since  we have $1+\varepsilon < \alpha <2$ then we get $1+\varepsilon < \frac{\alpha}{\alpha-1} < \frac{2}{\varepsilon}$.  Thus, on the one hand,  if the quantity $(a)$  above verifies $(a)<1$ then we have $\ds{\frac{1}{2}[(a)]^{\frac{2}{\varepsilon}} \leq  \frac{1}{2}[(a)]^{\frac{\alpha}{\alpha-1}}}$.  On the other hand, if the quantity $(a)$ verifies $(a)\geq 1$ then we have $\ds{\frac{1}{2}[(a)]^{1+\varepsilon} \leq  \frac{1}{2}[(a)]^{\frac{\alpha}{\alpha-1}}}$.

\end{appendices}

\section*{Data availability statement}
Data sharing not applicable to this article as no datasets were generated or analyzed during the current study.


\begin{thebibliography}{20}
\bibitem{BahouriDanchinCheman} H. Bahouri, J.Y. Chemin \& R. Danchin. \emph{Fourier Analysis and nonlinear partial differential equations}. Springer Vol: 343 (2011).
\bibitem{Biccari} U. Biccari \& V. Hern\'andez-Santama\'ia. \emph{The poisson equation from non-local to local}. Electronic Journal of Differential Equations, Vol. 2018 No. 145, pp. 1–13 (2018).
\bibitem{Biler1} P. Biler, G. Karch \& W. A. Woyczynski. \emph{Asymptotics for conservation laws involving Lévy diffusion generators}. Studia Math. 148, 171–192 (2001).	
\bibitem{BiTaWo} P. Biler, T. Funaki \& Wojbor A. Woyczynski. \emph{Fractal Burgers Equations}. 
Journal of Differential Equations,
	Volume 148, Issue 1,
Pages 9-46 (1998).   
\bibitem{BiKarWo} P. Biler, G. Karck \& W. Woyczy\'nski. \emph{Asymptotics for conservation laws involving Lévy diffusion generators}. Studia mathematica: 148 (2) (2001).
\bibitem{BrandKarch} L. Brandolese \& G. Karch. \emph{Far field asymptotics of solutions to convection equation with anomalous diffusion.} J. Evolution Equations. 8: 307–326 (2008).
\bibitem{Cui}  S. Cui \& X. Zhao. \emph{Well-posedness of the Cauchy problem for Ostrovsky, Stepanyams and Tsimring
equation with low regularity data}. J. Math. Anal. Appl. 344 778–787 (2008).
\bibitem{Dro}  J. Droniou, T. Gallou \&  J. Vovelle. \emph{Global solution and smoothing effect for a non-local regularization of a hyperbolic equation}, J. Evol. Eq. 3, 499–521  (2002).
\bibitem{Dro2} J. Droniou, C. Imbert. \emph{Fractal first order partial differential equations}. Arch. Rat. Mech. Anal. 182, 299–331 (2006).
\bibitem{FuSuWo}T. Funaki, D. Surgailis \& W. A. Woyczynski. \emph{Gibbs-Cox random fields and Burgers  turbulence}.  Ann. Appl. Prob. 5, 701-735 (1995).
\bibitem{FuWo2} T. Funaki \& W. A. Woyczynski. \emph{ Interacting particle approximation for fractal Burgers equation}.  Stochastic Processes and Related Topics, A Volume in Memory of Stamatis Cambanis, Birkha\"user, Boston (1998).
\bibitem{IgnatRossi} L. Ignat \& J. D. Rossi. \emph{A non-local convection-diffusion equation}. Journal of Functional Analysis, Volume 251, Issue 2:  399-437 
(2007).
\bibitem{Jacob} N. Jacob. \emph{Pseudo-differential operators and Markov processes}. Vol. I. Fourier analysis and semi-groups. Imperial College Press, London, (2001). 
\bibitem{Jourdan} B. Jourdain, S. Méléard \& W. A. Woyczynski. \emph{A probabilistic approach for nonlinear equations involving the fractional Laplacian and singular operator}.  Potential Analysis 23 , 55–81 (2005).
\bibitem{SaiWo1} A. S. Saichev \& W. A. Woyczynski. \emph{Advection of passive and reactive tracers in multi-dimensional Burgers velocity field}.  Physica D 100, 119-141 (1997).
\bibitem{SaiWo2} A. S. Saichev \& W. A. Woyczynski. \emph{Distributions in the Physical and Engineering Sciences}. Distributional and Fractal Calculus, Integral Transforms and Wavelets, Vol. 1, Birkha\"user, Boston, (1997).
\bibitem{ShZaFri}   M. F. Shlesinger, G. M. Zaslavsky, \& U. Frisch. \emph{L\'evy Flights and Related Topics
in Physics}. Lecture Notes in Physics, Vol. 450, Springer$-$Verlag, Berlin (1995). 
\bibitem{Vazquez}  J. L. V\'azquez, Arturo de Pablo, Fernando Quir\'os, Ana Rodr\'iguez. \emph{Classical solutions and higher regularity for nonlinear fractional diffusion equations}.    J. Eur. Math. Soc. 19 , 7: 1949-1975 (2017).
\bibitem{Zas1} G. M. Zaslavsky. \emph{Fractional kinetic equations for Hamiltonian chaos}.  Physica D 76, 110-122 (1994).
\bibitem{Zas2} G. M. Zaslavsky \& S. S. Abdullaev. \emph{Scaling properties and anomalous transport of particles inside the stochastic layer}. Phys. Rev. E 51, No. 5  3901-3910 (1995).
\end{thebibliography}
\end{document}